\documentclass{article}

\usepackage{amsfonts, bezier, amsmath,
amsxtra, amsbsy, amsgen, amsopn,
amssymb,
 verbatim, graphicx, mathrsfs, bm, bbm, epsfig}

\newtheorem{defn}{Definition}[section]

\newtheorem{lem}[defn]{Lemma}

\newtheorem{prop}[defn]{Proposition}
\newtheorem{thm}[defn]{Theorem}
\newtheorem{conj}[defn]{Conjecture}
\newtheorem{coro}[defn]{Corollary}
\newtheorem{constr}[defn]{Construction}
\newtheorem{definition}{Definition}[section]
\newtheorem{proof}[defn]{Proof}

\newtheorem{alg}{\noindent Algorith}[section]
\newenvironment{algo}[1][] {\begin{alg} [#1] $\;$} {\end{alg}}

\newtheorem{rem}[defn]{Remark}

\newtheorem{notation}{\noindent Notation}[section]

\newcommand{\R}{\textnormal{I}\mspace{-4mu}\textnormal{R}}

\newcommand{\bdefn}{\begin{defn}}
\newcommand{\edefn}{\end{defn}}
\newcommand{\blm}{\begin{lem}}
\newcommand{\elm}{\end{lem}}
\newcommand{\bprop}{\begin{prop}}
\newcommand{\eprop}{\end{prop}}
\newcommand{\bthm}{\begin{thm}}
\newcommand{\ethm}{\end{thm}}
\newcommand{\bconj}{\begin{conj}}
\newcommand{\econj}{\end{conj}}
\newcommand{\bcoro}{\begin{cor}}
\newcommand{\ecoro}{\end{cor}}
\newcommand{\balgo}{\begin{algo}}
\newcommand{\ealgo}{\end{algo}}
\newcommand{\brm}{\begin{rem}}
\newcommand{\erm}{\end{rem}}
\newcommand{\bnot}{\begin{notation}}
\newcommand{\enot}{\end{notation}}
\newcommand{\bpf}{\begin{proof}}
\newcommand{\epf}{\end{proof}}
\newcommand{\beg}{\begin{eg}}
\newcommand{\eeg}{\end{eg}}
\newcommand{\ben}{\begin{enumerate}\renewcommand{\theenumi}{(\roman{enumi})}}
\newcommand{\een}{\end{enumerate}}
\newcommand{\be}{\begin{eqnarray}}
\newcommand{\ee}{\end{eqnarray}}
\newcommand{\nn}{\nonumber}

\renewcommand{\choose}[2]{\genfrac{(}{)}{0pt}{}{#1}{#2}}

\numberwithin{equation}{section} \everymath{\displaystyle}

\newcommand{\beas}{\begin{eqnarray*}}
\newcommand{\enas}{\end{eqnarray*}}

\newcommand{\bea}{\begin{eqnarray}}
\newcommand{\ena}{\end{eqnarray}}

\newcommand{\bfw}{{\bf{W}}}

\newcommand{\bfb}{{\bf{b}}}

\begin{document}

\title{Joint Vertex Degrees in an Inhomogeneous Random Graph Model}
\author{K. Lin \footnote{Department of Statistics,
University of Oxford, 1 South Parks Road, OXFORD OX1 3TG, UK;   email: \texttt{k.lin@wolfson.oxon.org} 
\medskip}
\ and
G. Reinert\footnote{Department of Statistics,
University of Oxford, 1 South Parks Road, OXFORD OX1 3TG, UK;  email: \texttt{reinert@stats.ox.ac.uk}
} }

\date{} 
\maketitle


\begin{abstract}
In a random graph, counts for the number of vertices with given degrees will typically be dependent. We show via a multivariate normal and a Poisson process approximation that, for graphs which have independent edges, with a possibly inhomogeneous distribution, only when the degrees are large can we reasonably approximate the joint counts as independent. The proofs are based on Stein's method and the Stein-Chen method with a new size-biased coupling for such inhomogeneous random graphs, and hence bounds on distributional distance are obtained. Finally we  illustrate that apparent (pseudo-) power-law type behaviour can arise in such inhomogeneous networks despite not actually following a power-law degree distribution.
\end{abstract}

\section{Introduction}
It has been observed in many real-world networks that, when plotting the observed number of vertices of degree $k$ against $k$, on a log-log-scale the plots tend to look linear. This so-called {\it scale-free behaviour}, see e.g. \cite{dormenbook}, has motivated the {\it scale-free network model} introduced by \cite{barabasialbert}, yielding a probability distribution for the number of vertices of large degree which is scale-free.

Some issues arise when trying to assess the vertex degree distribution from a single network.  The log-log scale lends itself to over-interpretation; \cite{solowetal} discusses a good number of pitfalls when trying to test for power law using such plots. Moreover, \cite{stumpfetal} have shown that when sampling from a scale-free network, the sampled network will not in general be scale-free. In addition, the total number of vertices in the network is fixed, and hence counts for different degrees will be dependent. We shall see in this paper that whether or not the dependence is negligible depends on the size of the degrees under consideration - only when then degrees are large can we reasonably approximate the joint counts as independent. We establish these facts by proving a multivariate normal approximation, with possibly non-diagonal asymptotic covariance matrix,  as well as a Poisson process approximation, with independent coordinates. We give bounds for these approximations which depend on the size of the degrees under consideration. Finally we shall illustrate that apparent (pseudo-) power-law type behaviour can arise in networks which are constructed using independent edges, and do not follow a power-law behaviour.

The degree of a vertex is one of the fundamental summaries for random graphs, and hence the degree distribution is a natural object to study.  In a general  random graph $\mathscr{G}^n$ on a set $V$ of  $n$ vertices, the degree of a vertex $v$, denoted by $D(v)$, is defined as the number of vertices adjacent to $v$.
The most basic model of a random graph is that of  Bernoulli graph  $\mathscr{G}(n,p)$, introduced by Erd\"{o}s $\&$ R\'{e}nyi \cite{er}. A survey  of Poisson approximation for distribution of the $k$'th largest degree for large $k$ in the Bernoulli model $\mathscr{G}(n,p)$, as well as of both a Poisson approximation and a normal approximation for the number of vertices of a given degree can be found in \cite{bbb5}, with bounds on the distributional distance. For the joint distribution of degrees in the Bernoulli model $\mathscr{G}(n,p)$,  \cite{a} give an approximation with simpler models derived from a Binomial distribution and use this for univariate normal approximations.

While Bernoulli random graphs typically do not model real-world networks well, in \cite{daudinetal} a mixture model for Bernoulli random graphs is shown to be suitable for some biological networks. Under the name {\it stochastic block model} a similar mixture model has proven successful in the area of social network analysis, see \cite{nowickisnijders}. Here we use
the inhomogeneous model $\mathscr{G}(n,\{p_{ij}\})$  as a sub-model of $\mathscr{G}^n$, consisting of all graphs in which the edges occur independently, and for $i,j\in V$ the probability that vertices $i$ and $j$ are adjacent is $p_{ij}$. This general model not only includes Bernoulli random graphs, but also mixtures of Bernoulli random graphs, Newman-Moore-Watts-Strogatz small world networks as defined in \cite{nmw}, and the expontial random graph model, which is defined by assuming in $\mathscr{G}(n,\{p_{ij}\})$ that $p_{ij}=\exp{(\theta_i+\theta_j)}/\{1+\exp{(\theta_i+\theta_j)}\}$, where $\{\theta_i,i\in V\}$ are parameters of the model. For fairly general random graph models which include a
 Barab\'{a}si-Albert scale-free model, but do not quite cover the class $\mathscr{G}$\lowercase{$(n,\{p_{ij}\})}$ in full generality,  \cite{bollobasetal} give a univariate mixed Poisson approximation for the number of vertices with a given degree.
There is a lack of results for multivariate approximations, despite the need to understand log-log plots. In addition, networks consist of a finite number of vertices, and, depending on the complexity, the distribution of vertices with a fixed degree may be far from the asymptotic regime; thus bounds on the distributional approximations are required.

In order to understand log-log plots of the number of vertices with degree $k$ versus $k$, we consider the
 degree-count sequence $\mathbf{W}:=(W_i,0\leq i\leq n-1)$, where $W_i$ counts the number of vertices having degree exactly $i$. The definitions of both sequences $\mathbf{D} := (D(v), v \in V) $ and $\mathbf{W}$ can be related by introducing the index set
 \be
\Gamma:=\{(v,i):v\in V;0\leq i\leq n-1\},
\label{Gamma-all}
\ee
and defining, for $(v,i)\in\Gamma$, the Bernoulli random variables $X_{(v,i)}:=\mathbbm{1} ( D(v)=i )$, where $\mathbbm{1}( \cdot ) $ is the indicator function. Then $D(v)=\Sigma_{i=0}^{n-1} i X_{(v,i)}$, and  $W_i=\Sigma_{v\in V} X_{(v,i)}$. Other interesting statistics may be also obtained by this setting. For instance, one may define random variable $Z_k=\Sigma_{i\geq k}W_i$ as the number of vertices having degree at least $k$, for $0\leq k\leq n-1$, and consider the sequence $\mathbf{Z}:=(Z_k,0\leq k\leq n-1)$.

In the flavour of probability theory, as the sequences $\mathbf{D}$, $\mathbf{W}$ and $\mathbf{Z}$ are deterministic functions of the collection $\mathbf{X}:=\{X_{(v,i)},(v,i)\in\Gamma\}$, the $\sigma$-fields $\sigma(\mathbf{D})$, $\sigma(\mathbf{W})$ and $\sigma(\mathbf{Z})$, generated by $\mathbf{D}$, $\mathbf{W}$ and $\mathbf{Z}$ respectively, are all contained in the $\sigma$-field $\sigma(\mathbf{X})$ generated by $\mathbf{X}$. The collection $\mathbf{X}$ in turn can be represented by the point process $\Xi$ defined by
\beas
\Xi:=\sum_{\alpha\in\Gamma}\delta_{\alpha}X_{\alpha},
\enas
where $\delta_{\alpha}$ is the point measure at $\alpha$, that is, for a set $B$, $\delta_{\alpha}(B)=1$ if $\alpha\in B$, or otherwise $\delta_{\alpha}(B)=0$.

For the degree-count sequence $\mathbf{W}=(W_i,0\leq i\leq n-1)$ in $\mathscr{G}(n,p)$,  two results in particular have inspired the current work. \cite{bhj} give  univariate Poisson approximations for the distribution of $W_k$ and $Z_k$, and \cite{GoldsteinRinott1996} prove a multivariate normal approximation for  the joint distribution of any sub-sequence $(W_{d_1},W_{d_2},\ldots, W_{d_m})$ of $\mathbf{W}$. Both results use Stein's method; in the context of Poisson approximation this method is usually called the {\it{Stein-Chen method}}. The applications of Stein's method in these two papers use a coupling construction to compute bounds on the errors made in the distributional approximations. For graph degrees counts, for
any $\alpha \in \Gamma$  a new graph model $\mathscr{G}^{\alpha}(n,p)$ is constructed, conditional on the model $\mathscr{G}(n,p)$, such that the distribution of $\mathscr{G}^{\alpha}(n,p)$ is the same as the conditional distribution of $\mathscr{G}(n,p)$ given $X_\alpha=1$; this coupling is a special case of a {\it size-bias coupling}. The difference between the degree-counts in $\mathscr{G}(n,p)$ and in $\mathscr{G}^{\alpha}(n,p)$ is then used ingeniously to give a bound on the distance to the target distribution.

In Section \ref{coupling} we construct such a coupling in the inhomogeneous model $\mathscr{G}(n,\{p_{ij}\})$, generalizing the existing construction for the homogeneous model. This coupling will be the main tool for our distributional approximations, which we derive in Section \ref{approximations}.
Firstly, in Theorem \ref{mvn}, we provide a multivariate normal approximation for the joint counts of vertices with pre-described degrees. The bound depends on the chosen degrees, and on the heterogeneity of the underlying graph. The approximating normal distribution has non-diagonal covariance matrix in general, and hence in the normal limit the counts will often  not be independent.

The multivariate normal approximation is suitable when the degrees under consideration are not too far away from the centre of the degree distribution. For large degrees, a compound Poisson approximation is more appropriate. Indeed Theorem \ref{PoissonThm} gives a Poisson point process approximation for the {\it{$M$-truncated}} point process $\Xi_M$ defined by
\be \label{Xi-M}
\Xi_M:=\sum_{\alpha\in\Gamma_M}\delta_\alpha X_\alpha,
\ee
where  for $0\leq M\leq n-1$, we put $\Gamma_M:=\{(v,i):v\in V;M\leq i\leq n-1\}.$

Using the invariant property of the total variation distance in functional transformations of point processes, we obtain, from Theorem \ref{PoissonThm}, in Corollary \ref{PoissonCoro1}  a multivariate compound Poisson approximation for the $M$-truncated degree sequence $\mathbf{D}_M := (D(v) \mathbbm{1} (D(v) \ge M), v \in V)$ in $\mathscr{G}(n,\{p_{ij}\})$. The result shows that counts for large vertex degrees are asymptotically independent when the edge probabilities are not too heterogeneous. All these results also contain a bound in distributional distance. This bound depends on the size of the degrees under consideration, and on the number of vertices, as well as on the heterogeneity in the edge probabilities.

We illustrate our results using simulations for a Bernoulli random graph as well as several classes of inhomogeneous random graphs.
Finally we show that the log-log plots for vertex degrees can appear to be power-law like, without following a power law, when the edge probabilities are small.

Proofs are postponed until Section \ref{proofs}.

\section{A Size Biased Coupling for Vertex Degrees} \label{coupling}

The size-biased distribution of a collection of variables $\mathbf{X}$  relates to a sampling procedure where the probability of an item to be included in the sample is directly proportional to its size. Formally it can be defined as follows, see
 for example, \cite{GoldsteinRinott1996}.

\begin{definition} Let $\mathcal{I}$ be an arbitrary index set and let $\mathbf{X} = \{ X_\alpha: \alpha \in \mathcal{I} \}$ be a collection of non-negative random variables with  means $\mathbb{E} X_\alpha = \lambda_\alpha > 0$. For $\beta \in \mathcal{I}$, we say that $\mathbf{X}^\beta = \{ X_\alpha^\beta: \alpha \in \mathcal{I} \}$ has \emph{the $\mathbf{X}$-size biased distribution in the $\beta^{th}$ coordinate} if
\be
\mathbb{E} X_\beta G(\mathbf{X}) = \lambda_\beta \mathbb{E} G(\mathbf{X}^\beta)  \nn
\ee
for all functions $G$ such that the expectations exist.
\label{defn-sizebias}
\end{definition}


A construction of  $(\mathbf{X}, \mathbf{X}^\beta)$, for each $\beta \in \mathcal{I}$, on a joint probability space is called a {\it{size-biased coupling}}.
 For any subset $B \subset \mathcal{I}$, we set $X_B = \Sigma_{\alpha \in B} X_\alpha$, and $\lambda_B = \mathbb{E} X_B$.
\cite{GoldsteinRinott1996} give the following mixture construction of a size-biased coupling for $\mathbf{X}$ in ``coordinate'' $B$:
Suppose  that  $\lambda_B < \infty$, and that for $\beta \in B$, we have a variable $\mathbf{X}^\beta$ which has the $\mathbf{X}$-size biased distribution in coordinate $\beta$ as in Definition \ref{defn-sizebias}. Then the random variable $\mathbf{X}^B$ which is obtained as the mixture of the distributions $\mathbf{X}^\beta$, $\beta \in B$ with weights $\lambda_\beta / \lambda_B$, satisfies that $
\mathbb{E} X_B G(\mathbf{X}) = \lambda_B \mathbb{E} G(\mathbf{X}^B)$.

The application of this construction for coupling variables for the degree-count sequence $\mathbf{W} = (W_{d_i}, 1\leq i\leq m)$ has been carried out in \cite{GoldsteinRinott1996}, for
distinct and fixed $d_i$ (where $i=1,\ldots,m$). The idea is, for a given degree $d_i$ and vertex $v$, we force $v$ to have degree $d_i$. If the degree of $v$ was equal to $d_i$ in the first place, no adjustment is necessary. If the degree $D(v)$ of $v$ in the original graph was larger than $d_i$, then $D(v) - d_i$ edges are chosen at random from the edges which include $v$ as one end point, and are removed. If $D(v) < d_i$, then $d_i  - D(v)$ edges of the form $\{u,v\}$ are added to the graph, where the vertices $u$ are chosen uniformly at random from the $n-1-D(v)$ vertices not adjacent to $v$.  Randomizing over the pair $(v, d_i)$, for $v\in V$ and $i=1,\ldots,m$, then gives a size-biased version of the graph $\mathscr{G}(n,p)$.


For the inhomogeneous model $\mathscr{G}(n, \{ p_{ij} \})$ we use the index set $\Gamma$ in (\ref{Gamma-all}), which covers all possible combinations between vertices and their degrees, and we write $A_i = V \times \{i\}$ for $i \in \{0, \ldots, n-1\}$. For $i \in \{0, \ldots, n-1\}$ we construct $(\mathbf{X}, \mathbf{X}^\beta)$ for $\beta \in A_i$, with $\mathbf{X} = \{X_{(v,i)}: (v,i) \in \Gamma\}$ and $X_{(v,i)} = \mathbbm{1}( D(v)=i )$, as detailed below; we call the resulting graph $\mathscr{G}^{\beta}(n, \{ p_{ij} \})$.

Then, we construct $(\mathbf{W}, \mathbf{W}^i)$ by using a random index $I$ in $A_i$, which has the probability mass function $\mathbb{P}(I=\beta) = \lambda_\beta / \lambda_{A_i}$, independently of all other random variables in the system; we call the resulting graph $\mathscr{G}^{i}(n, \{ p_{ij} \})$.

To describe the detailed construction of $\mathscr{G}^{\beta}(n, \{ p_{ij} \})$, let $E$ denote the (potential) edge set of the graph model, and define, for an edge $\{a,b\}\in E$, the Bernoulli random variables $X_{\{a,b\}}:=\mathbbm{1}(a\sim b)$ in $\mathscr{G}(n, \{ p_{ij} \})$ (and similarly $X_{\{a,b\}}^{\beta}$ in $\mathscr{G}^{\beta}(n, \{ p_{ij} \})$), where $a \sim b$ denotes the event that $a$ is adjacent to $b$. We also use the following notation: $N(v)$ the random neighbourhood of vertex $v$ in $\mathscr{G}(n, \{ p_{ij} \})$ (similarly $N^{\beta}(v)$ in $\mathscr{G}^{\beta}(n, \{ p_{ij} \})$), and $\mathbf{x}_i=\mathbf{x}_i(v)$ a $i$-set (i.e. a set with $i$ elements) of $V_v:=V\setminus\{v\}$.

\begin{constr}  \label{II}
For each $\beta=(v,i)\in A_i$, conditional on $\mathscr{G}(n, \{ p_{ij} \})$:

If $D(v)=d=i$, let $\mathscr{G}^{(v,i)}(n,\{p_{ij}\})=\mathscr{G}(n,\{p_{ij}\})$, that is,
 set $X_{\{a,b\}}^{(v,i)}=X_{\{a,b\}}$ for all $\{a,b\}\in E$.

If $D(v)=d>i$ and $N(v)=\mathbf{x}_d$, then we choose $\mathbf{x}_i\subset\mathbf{x}_d$ with probability proportional to
\be
\lefteqn{f^{+}(\mathbf{x}_i\,|\,\mathbf{x}_d)} \label{delp} \\
&:=&\sum_{j=0}^{i}\frac{1}{\choose{d-j}{i-j}}
\mathbb{P}\big(\,|N(v)\cap\mathbf{x}_i|=j,\;|N(v)\cap(V_v\setminus\mathbf{x}_d)|=i-j\;\big|\;D(v)=i\,\big) \nn
\ee
and delete all the edges between $v$ and the vertices in $\mathbf{x}_d\setminus\mathbf{x}_i$.
That is, with probability proportional to (\ref{delp}), we set $X_{\{v,x\}}^{(v,i)}=0$ for $x$ in $\mathbf{x}_d\setminus\mathbf{x}_i$, and $X_{\{a,b\}}^{(v,i)}=X_{\{a,b\}}$ for all $\{a,b\}$ elsewhere.

If $D(v)=d<i$ and $N(v)=\mathbf{x}_d$, then we choose $\mathbf{x}_i\supset\mathbf{x}_d$ with probability proportional to
\be
\lefteqn{f^{-}(\mathbf{x}_i\,|\,\mathbf{x}_d)} \label{addp} \\
&:=&\sum_{j=i}^{n-1}\frac{1}{\choose{j-d}{j-i}} \nn \\
&&\times 
\mathbb{P}\big(\,|N(v)\cup\mathbf{x}_i|=j,\;|N(v)\cup(V_v\setminus\mathbf{x}_d)|=n-1-j+i\;\big|\;D(v)=i\,\big), \nn
\ee
and add all the edges between $v$ and the vertices in $\mathbf{x}_i\setminus\mathbf{x}_d$.
That is, with probability proportional to (\ref{addp}), we set $X_{\{v,x\}}^{(v,i)}=1$ for $x$ in $\mathbf{x}_i\setminus\mathbf{x}_d$, and $X_{\{a,b\}}^{(v,i)}=X_{\{a,b\}}$ for all $\{a,b\}$ elsewhere.
\end{constr}

In Section \ref{proofs} we shall prove that \eqref{delp} and \eqref{addp} indeed are probabilities.
\begin{lem} \label{isaprob}
We have that
\be
\sum_{\mathbf{x}_i:\mathbf{x}_i\subset\mathbf{x}_d}f^{+}(\mathbf{x}_i\,|\,\mathbf{x}_d)
=1 \mbox{ and }
\sum_{\mathbf{x}_i:\mathbf{x}_i\supset\mathbf{x}_d}f^{-}(\mathbf{x}_i\,|\,\mathbf{x}_d)
=1.  \nn
\ee
\end{lem}

\brm
Note that, in all cases, the above construction \ref{II}  yields indeed that $N^{(v,i)}(v)=\mathbf{x}_i$ if $N(v)=\mathbf{x}_d$.  
\erm

We shall show in Section \ref{proofs}  that the distribution of $\mathscr{G}^{\beta}(n,\{p_{ij}\})$ is indeed the same as the conditional distribution of $\mathscr{G}(n,\{p_{ij}\})$ given $X_\beta=1$, yielding a construction of $(\mathbf{X}, \mathbf{X}^\beta)$ for $\beta\in A_i$, which in turn gives a  construction of $(\mathbf{W}, \mathbf{W}^i)$ via $\mathscr{G}^{i}(n, \{ p_{ij} \})$ using the random index $I\in A_i$ with $\mathbb{P}(I=\beta) = \lambda_\beta / \lambda_{A_i}$. In the next section, we shall use Construction \ref{II} in $\mathscr{G}(n,\{p_{ij}\})$ to obtain a multivariate normal approximation for the degree-count sequence $\mathbf{W}$, and a compound Poisson approximation for the truncated degree sequence $\mathbf{D}_M$.

\section{Approximations for Degree Counts}

\subsection{Multivariate Normal Approximation}

For a multivariate normal approximation we generalize the argument from \cite{GoldsteinRinott1996}, which is based on Stein's method. Let $V=\{1,\ldots,n\}$,  let $d_i, i=1, \ldots, p,$ be distinct numbers
in $\{0, 1, \ldots, n-1\}$, and let $W_{d_i} = \Sigma_{v=1}^n \mathbbm{1}(D(v) = d_i) $ denote the number of vertices of degree $d_i$ in
$\mathscr{G}(n, \{ p_{ij} \})$. Denote by
${\bf W} = (W_{d_1}, \ldots, W_{d_p})$ the vector of degree counts. As $D(v)$ has a Poisson-binomial distribution which is cumbersome to
write explicitly, we abbreviate
\beas
q_{v, d} := \mathbb{P}(D(v) = d).\nn
\enas
Let  $\bm{\lambda} = (\lambda_1, \ldots, \lambda_p)$ denote the expectation vector of $\bfw$; for $i=1, \ldots, p$,
\beas
\lambda_i := \mathbb{E}W_{d_i} = \sum_{v=1}^n \mathbb{P}(D(v) = d_i) = \sum_{v=1}^n q_{v, d_i}.\nn
\enas
We also abbreviate, for $v = 1, \ldots, n$,
\beas 
\mu_v := \mathbb{E}D(v) = \sum_{u} p_{u,v},
\enas
where, and as everywhere else, we use the convention $p_{v,v}=0$ in $\mathscr{G}(n, \{ p_{ij} \})$. Then
\bea \label{vard}
Var(D(v)) & =& \sum_{u} p_{u,v}(1 - p_{u,v}) \leq \mu_v.
\ena

For vertices $v_1,v_2,\ldots,v_m$, denote by $\mathscr{G}^{(v_1,\ldots,v_m)}(n-m, \{ p_{ij} \})$ the random graph $\mathscr{G}(n, \{ p_{ij} \})$ with vertices $v_1,v_2,\ldots,v_m$ and all their edges removed. For this graph let $D^{(v_1,\ldots,v_m)}(w)$ be the degree of vertex $w$, where $w\notin\{v_1,\ldots,v_m\}$. We let
\beas 
q_{w, d}^{(v_1,\ldots,v_m)} := \mathbb{P}(D^{(v_1,\ldots,v_m)}(w) = d).
\enas
It is straightforward to calculate that the entries of $ \Sigma= (\sigma_{i,j})$, the covariance matrix of ${\bf W}$, are
\beas
\sigma_{ij}
&=&  \mathbbm{1}(i=j) \lambda_i - \sum_{v}q_{v,d_i}q_{v,d_j} \\
&&+ \sum_v \sum_{w \ne v} \Big[ p_{wv} q_{v, d_i-1}^{(w)} q_{w, d_j-1}^{(v)} + (1-p_{wv}) q_{v, d_i}^{(w)} q_{w, d_j}^{(v)} - q_{v,d_i}q_{w,d_j} \Big].\nn
\enas

When the degrees are``typical", a multivariate normal approximation for ${\bf W} $ holds, as the following theorem shows.
Using
the notation from \cite{GoldsteinRinott1996},
for smooth functions $h: \R^p \rightarrow \R$,
we let $Dh$ denote the vector of first partial derivatives of $h$, and in general $D^k$ the $k^{th}$ derivative of $h$; $\parallel h \parallel$ denotes the supremum norm. We also abbreviate
\beas
{\bar{p}}_v = \frac{1}{n-1} \sum_{w} p_{wv}.
\enas
 Define $\Sigma_0=(\sigma_{ij}^0)$ with $\sigma_{ij}^0$ given by
\beas
\sigma_{ij}^0 &=& \mathbbm{1}(i=j) \lambda_i - \sum_{v}q_{v,d_i}q_{v,d_j} \\
&& 
+ \sum_{v, w} \sqrt{{\bar{p}}_v} \sqrt{{\bar{p}}_w}(q_{v,d_i-1}-q_{v,d_i})(q_{w,d_j-1}-q_{w,d_j}).
\enas
Finally let
\beas
M= \max\Big\{ \sum_v \mu_v; \sum_v \mu_v^3\Big\}.
\enas
The following result gives a bound for the distance between the distribution of our degree count vector $\bfw$ to a multivariate normal distribution with the same mean as $\bfw$, but with covariance matrix $\Sigma_0$.
The proof can be found in Section \ref{proofs}.

\begin{thm}\label{mvn}
For any function $h: \R^p \rightarrow \R$ having bounded mixed partial derivatives up to order 3,
\beas
\lefteqn{ \left| \mathbb{E} h( \Sigma_0^{-1/2} (\bfw - \bm{\lambda}) ) - Nh   \right| } \nn\\
&\le
p^3 \tau^2  \parallel D^2 h \parallel
\left( \sum_{i=1}^p  B_i  + S\right)  + \frac{p^5}{3}
\tau^{3}
\parallel D^3 h \parallel \left\{ M +  \sum_{i=1}^p   \lambda_i  (d_i + 1)^2  \right\}.\nn
\enas
Here
\beas \label{B1}
B_{i}&=& 128 \left\{ (10 + 6d_i^2) M + (d_i^2 + 2) \sum_{u,v} p_{u,v}^2 (3 M +1) \right\}^\frac12
\enas
and
\beas 
S&=&  4 \sum_{v,w} p_{vw}^2+  \sum_v \sum_{w\ne v  }  | p_{wv} - \sqrt{{\bar{p}}_v} \sqrt{{\bar{p}}_w}| | (q_{v,d_i-1}-q_{v,d_i})(q_{w,d_j-1}-q_{w,d_j}) | \nonumber \\
&& + \sum_v {\bar{p}}_v |  q_{v,d_i-1}-q_{v,d_i} | | q_{v,d_j-1}-q_{v,d_j}|
\enas
as well as
\beas 
\tau&=& {\Big[\sum_v \min_{i} q_{v,d_i}(1-\sum_{i}q_{v,d_i})\Big]^{-1/2}} .
\enas
\end{thm}

\begin{rem}
For every $v=1, \ldots, n$, the degree  $D(v)$ can be approximated by a Poisson distribution with
parameter $\mu_v$. From \cite{bhj}, Equation (1.23), p. 8,
\be \label{bhjbound}
d_{TV}\big( {\cal L}(D(v)), \mbox{Po}(\mu_v)\big) \leq \min\left(1; \frac{1}{\mu_v}\right) \sum_{u} p_{u,v}^2.
\ee
Here $d_{TV}$ denotes the total variation distance; for two probability measures $\mu$ and $\nu$ on the same probability space with $\sigma-$algebra ${\mathcal{B}}$, we define
\beas
d_{TV}(\mu, \nu) = \sup_{B \in {\mathcal{B}} } | \mu(B) - \nu(B)|.
\enas
The Poisson approximation is good for example when $p_{u,v} \approx \frac{c}{n-1} =\pi$ for all $u,v$ and some constant $c$;  then
$\mu_v \approx c = O(1)$, and   $\Sigma_{u} p_{u,v}^2 \approx \frac{c^2}{n-1}$.
The distributional regime where the normal approximation is plausible is when all degrees are moderate, $\mu_v = O(1)$
for all $v$, so that $M = O(n)$;   then it is reasonable to think of $\lambda_i \asymp n$ and  $q_{v, d_i-1}^{(u)}= O(1) $ as well as $ q_{v, d_i} = O(1)$.  In this regime, with $p$ fixed,
$
\sum_i B_{i}=  O\left( n^{\frac12} \left( 1 + \sum_i d_i \right)  \right) ,
$
and $\tau \asymp n^{-\frac12} $. If $\sum_i d_i^2  =O(1)$,  this yields an overall bound of the order
$n^{-\frac12}$.
\end{rem}

\begin{rem}
The term $S$ arises from the variance approximation;
\beas
  \sum_v \sum_{w }  | p_{wv} - \sqrt{{\bar{p}}_v} \sqrt{{\bar{p}}_w}| | (  q_{v,d_i-1}-q_{v,d_i})(q_{w,d_j-1}-q_{w,d_j}) |
 \enas
 vanishes when all $p_{wv} = \pi$ are equal.
\end{rem}

\begin{rem}
In the case that $p_{u,v} = \frac{c}{n}$ for all $u \ne v$, now putting $q_d =q_{v,d}$,  the approximating covariance simplifies to
\beas
\sigma_{ij}^0 &=& \mathbbm{1}(i=j) n q_{d_i} + n q_{d_i} q_{d_j} \left[ \frac{(n-1)(d_i-c)(d_j-c)}{n c \left( 1 - \frac{d_i}{n}\right) \left( 1 - \frac{d_j}{n}\right)  } -1 \right] .
\enas
Under the regime that $d_i$ and $d_j$ are typical degrees, so that $q_{d_i} $ and $q_{d_i} $ are moderate, this
expression will not in general tend to zero for $i \ne j$ as $n \rightarrow \infty$; the covariance does in general not vanish, and the degree counts will be asymptotically dependent.

\end{rem}

\begin{rem}
While our bounds are for smooth test functions $h$ only, they could be generalised to non-smooth test functions along the lines of \cite{Rinott1996}. Corresponding work is in progress, \cite{larryjay}, and to avoid duplicate work we restrict ourselves to smooth test functions.
\end{rem}

\subsection{Poisson Process Approximation for the Truncated Degree Sequences in $\mathscr{G}(n, \{ p_{ij} \})$}

The construction \ref{II} of $(\mathbf{X}, \mathbf{X}^\beta)$ allows to assess the distribution of the $M$-truncated degree sequence $\mathbf{D}_M:=(D(v)\mathbbm{1}(D(v)\geq M),v\in V)$ for an arbitrary integer $M\geq 0$. To this end, define $\Gamma_M$ by $\{1,\dots,n\}\times\{M,\ldots,n-1\}$, a subset of $\Gamma$, and restrict the definitions of $\mathbf{X}$ and $\mathbf{X}^\beta$ to $\Gamma_M$ to have $\mathbf{X}_M:=\{X_{(v,i)}:(v,i)\in\Gamma_M\}$ and $\mathbf{X}^\beta_M:=\{X_{(v,i)}^\beta:(v,i)\in\Gamma_M\}$, such that
$
\mathcal{L}(\mathbf{X}^\beta_M) = \mathcal{L}(\mathbf{X}_M \,|\, X_\beta=1),
$
for $\beta\in\Gamma_M$. Construction \ref{II}  can be used to derive a Poisson process approximation, with respect to the total variation distance, for the point process $\Xi_M$ defined in (\ref{Xi-M}), where the target Poisson point process $\Theta_M$ on $\Gamma_M$ has intensity $\bm{\lambda}_M=(\mathbb{E}X_\alpha,\alpha\in\Gamma_M)$.

\bthm
In $\mathscr{G}(n,\{p_{ij}\})$, we have
\be d_{TV}\big(\mathcal{L}(\Xi_M), \mbox{\emph{Po}}(\bm{\lambda}_M)\big)\leq b_{M,1}+b_{M,2}, \nn\ee
where
\be
b_{M,1}& =& \sum_{v\in V}\left( \mathbb{P}(D(v) \ge M )\right)^2, \mbox{ and } \nn \\
b_{M,2}& = &  2\sum_{v\in V}\sum_{u\in V_v}
\mathbb{P}(D(v)\geq M)\,\mathbb{P}(D^{(v)}(u)\geq M -1).  \nn
\ee
\label{PoissonThm}
\ethm

Since the total variation distance between the two processes also serves as an upper bound of the total variation distance between deterministic functions of the two processes, that is,
$d_{TV}\big(\mathcal{L}(f(\Xi_M)), \mathcal{L}(f(\Theta_M))\big)\leq d_{TV}\big(\mathcal{L}(\Xi_M), \mathcal{L}(\Theta_M)\big),$
where $\Theta_M\sim\mbox{Po}(\bm{\lambda}_M)$ and $f$ is any deterministic function, we assess the distribution of $\mathbf{D}_M$ by taking the function $f$ on point measures $\xi$ on $\Gamma_M$ as $f(\xi):=\Big(\sum_{i=0}^{n-1}i\xi((v,i)),v\in V\Big)$.
In this way, the target distribution ${\mathcal {L}} ( f(\Theta_M) ) $ gives rise to a multivariate compound Poisson approximation for $\mathbf{D}_M$, in the next corollary.
The result justifies  the independence assumption among large vertex degrees as used when interpreting log-log plots for vertex degrees when the degrees are observed not in independent graphs, but in the same graph. It also bounds the departure from an independent point process in terms of the degree threshold $M$.

\begin{coro} In the model $\mathscr{G}(n,\{p_{ij}\})$, let $\mathbf{D}_M$ denote the $M$-truncated degree sequence $(D(v)\mathbbm{1}(D(v)\geq M),v\in V)$, and $\mathbf{Y}_M$ denote the compound Poisson vector $(Y_{v,M},v\in V)$, in which all components are independent and $Y_{v,M}=\Sigma_{i=M}^{n-1}iY_{v,i}$ with $Y_{v,i}\sim\mbox{\emph{Po}}(\mathbb{E}X_{(v,i)})$. Then, with the $b_{M,i}$'s from Theorem \ref{PoissonThm},
\beas
d_{TV}\big(\mathcal{L}(\mathbf{D}_M),\mathcal{L}(\mathbf{Y}_M)\big)\leq b_{M,1}+b_{M,2}.
\enas
\label{PoissonCoro1}
\end{coro}

\begin{rem} Corollary \ref{PoissonCoro1} is consistent with Theorem 3.13 given in Bollob\'{a}s, Janson $\&$ Riordan (2007), where, in a fairly general sub-model of $\mathscr{G}(n,\{p_{ij}\})$, it is shown that the distribution of an individual vertex degree converges to a mixed Poisson distribution. In contrast, Corollary \ref{PoissonCoro1} not only applies in the multivariate case, but it also provides an explicit error bound on the distance.
\end{rem} 

\brm Using  $\mu_{u}^{(v)}=\Sigma_{x\in V\setminus\{u,v\}}p_{ux}$, we obtain from \eqref{bhjbound} that
\be
\mathbb{P}(D(v)\geq M) &\leq& \mbox{Po}(\mu_v)\{[M,n-1]\}+\frac{1-e^{-\mu_v}}{2\mu_v}\sum_{x\in V_v}p_{vx}^2, \mbox{ and } \nn \\
\mathbb{P}(D^{(v)}(u)\geq M-1) &\leq& \mbox{Po}(\mu_{u}^{(v)})\{[M-1,n-2]\}+\frac{1-e^{-\mu_{u}^{(v)}}}{2\mu_{u}^{(v)}}\sum_{x\in V\setminus\{u,v\}}p_{ux}^2.\nn
\ee
This yields an upper bound for the quantities in Theorem \ref{PoissonThm}; we can use the Poisson distribution as a guideline for a good choice of $M$. These probabilities could be further bounded using Proposition A.2.3 in \cite{bhj}.
\erm

\subsection{Simulations for the Correlation between Counts}

We now illustrate the dependence structure in four different random graph models, all on $n=100$ vertices, with independent edges. We estimate the correlations from  10,000 samples of graphs for each model. The models are as follows.

\begin{enumerate}
\item  {\bf M1}. The first model is the Bernoulli random graph with $p_{u,v} = p = \frac{1}{n}$. This graph is at criticality; some, but not all, realizations may yield a giant component, see \cite{bbb5}.
\item   {\bf M2}. In this model, $p_{u,v} = \frac{1}{5}$ if $0 < | u - v |\,(mod\;100) \le 10$, and $p_{u,v} = \frac{1}{80}$ if $ | u - v |\, (mod\;100) > 10$. This is a modified Newman-Moore-Watts small-world model, see \cite{nmw}, with $100$ vertices; two vertices  at most distance $10$ away from each other are connected with probability $\frac{1}{5}$, and two vertices more than distance $k$ away from each other are connected with probability $\frac{1}{80}$.

\item  {\bf M3}. Here $p_{u,v}= \frac{\min(u,v)}{n}$ for $u \ne v$; the smaller of the two vertices determines the probability.
\item  {\bf M4}. This model is motivated by Rasch-type models;
for $u<v$, we set $p_{uv}=\alpha_u(3)\alpha_v(10)$, with
\be
\alpha_u(i)&=&\Bigg\{ \begin{array}{l} 1/(i\sqrt{n}), \quad u\leq \frac{n}{2} \\ i/\sqrt{n},\qquad\; u> \frac{n}{2}  \end{array} \nn .
\ee
\end{enumerate}

\begin{figure}
\begin{center}
\includegraphics[height=5cm]{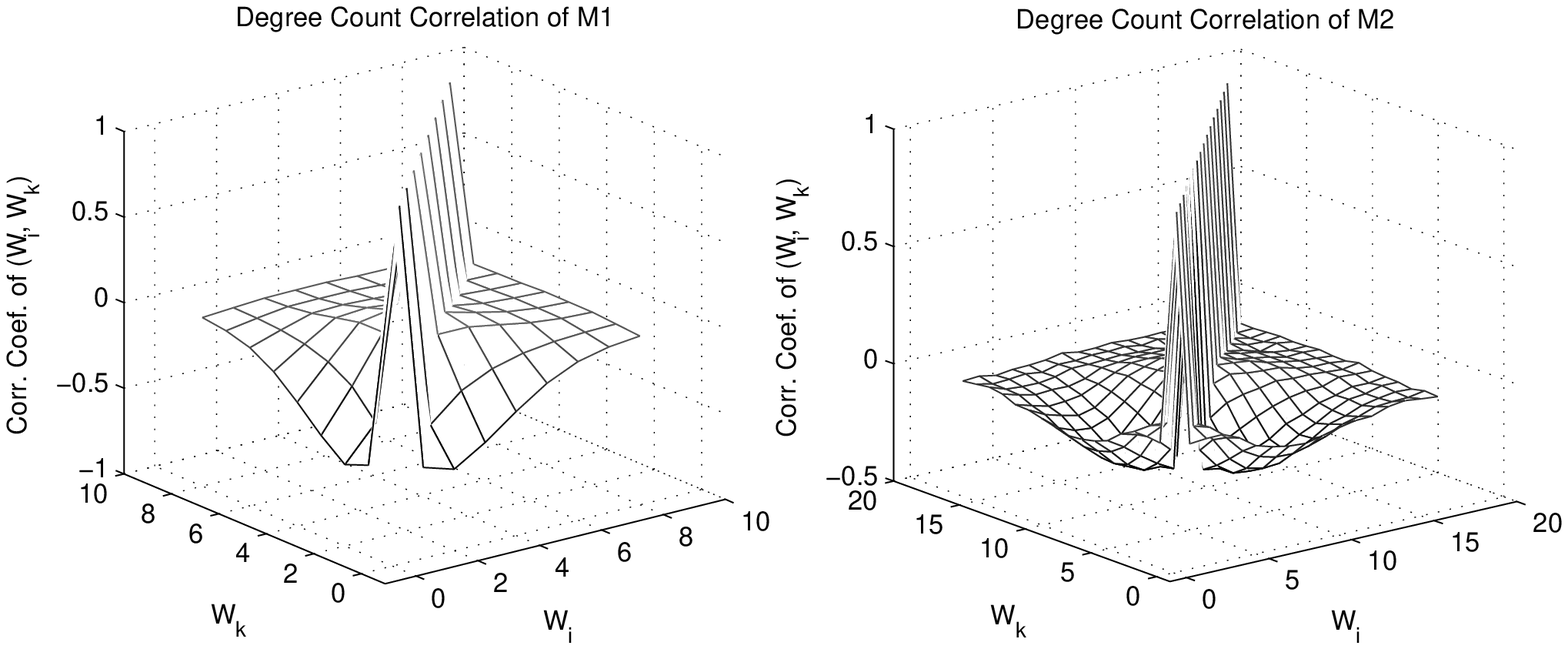}
\includegraphics[height=5cm]{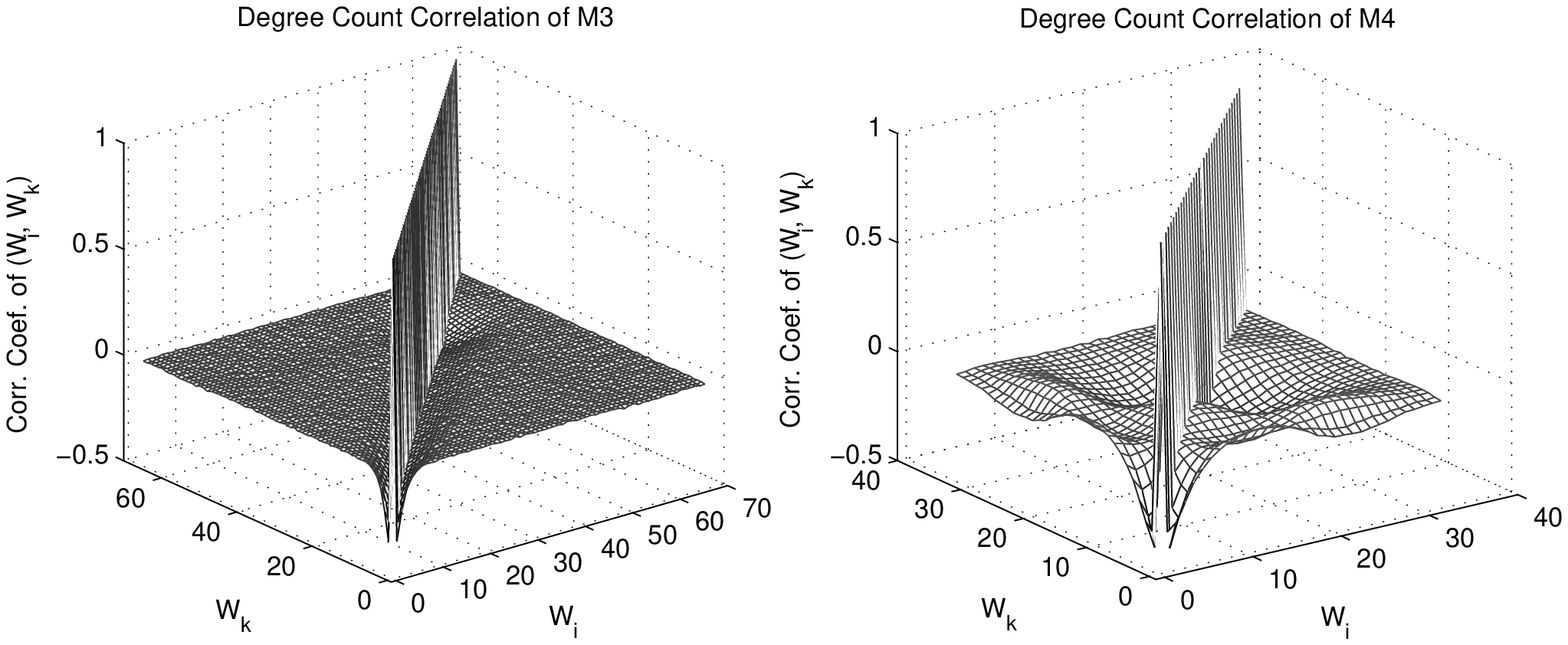}
\caption{\emph{Degree count correlation  in Models M1 -- M4. }} \label{dccorr}
\end{center}
\end{figure}

Figure \ref{dccorr} shows the correlations between the degree counts in the four models; except for Model M3  there is an appreciable correlation even far away from the diagonal.

\begin{figure}
\begin{center}
\includegraphics[height=4cm]{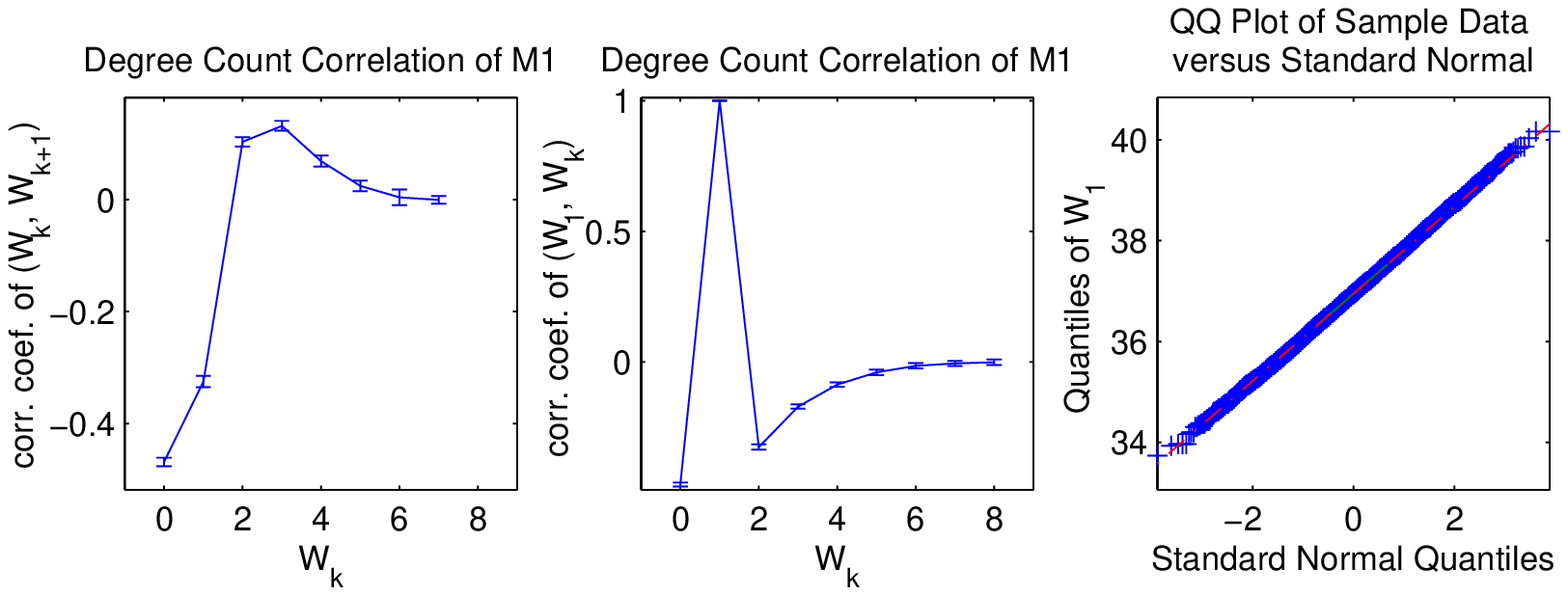}
\includegraphics[height=4cm]{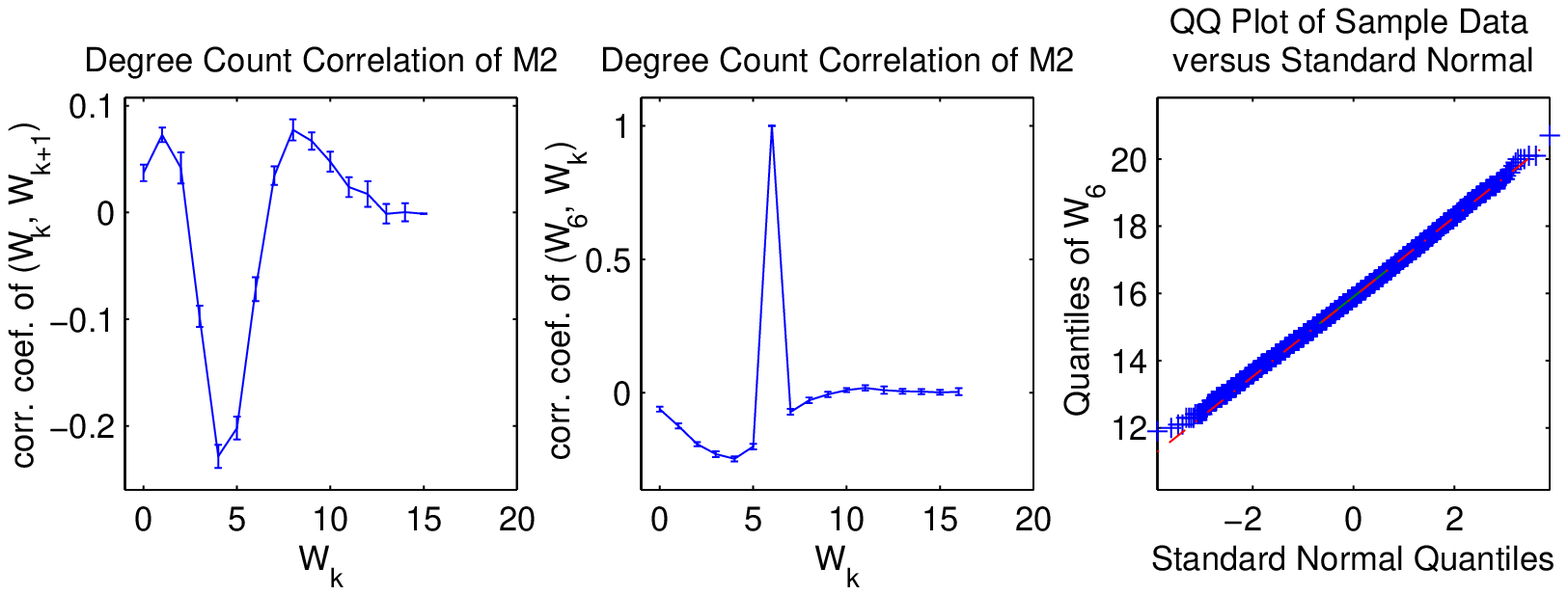}
\includegraphics[height=4cm]{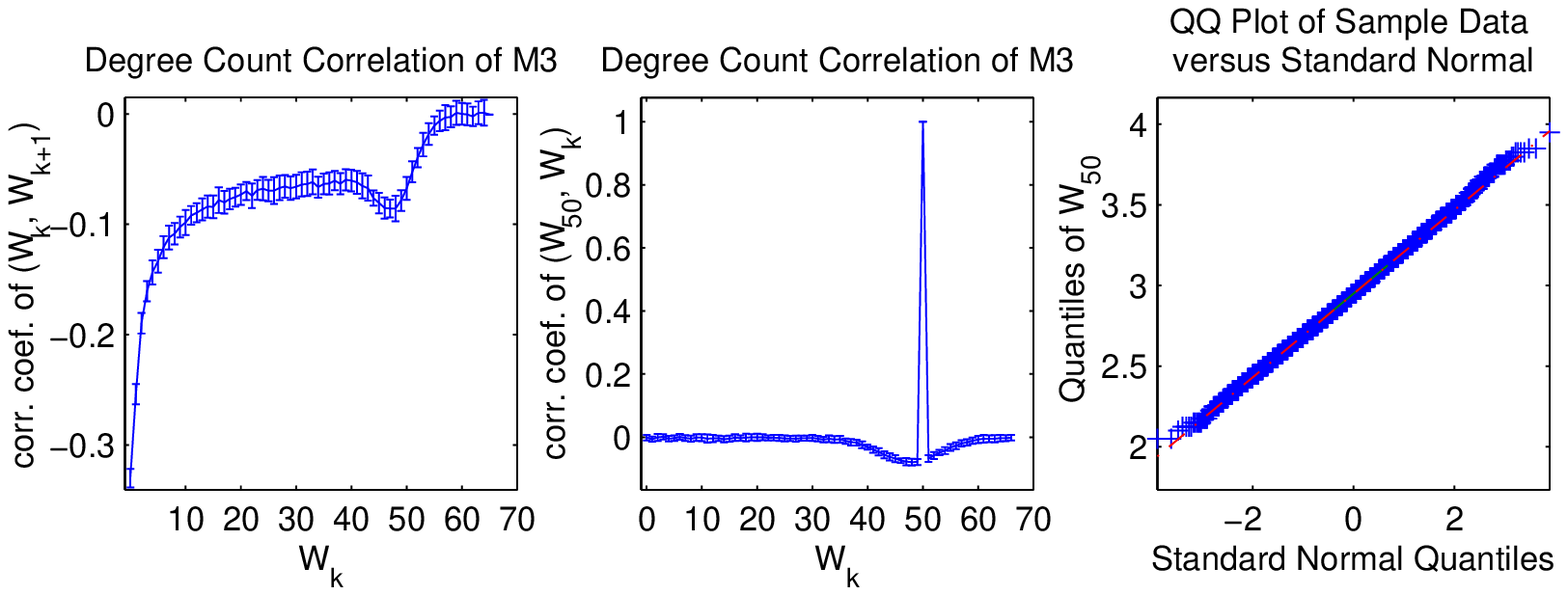}
\includegraphics[height=4cm]{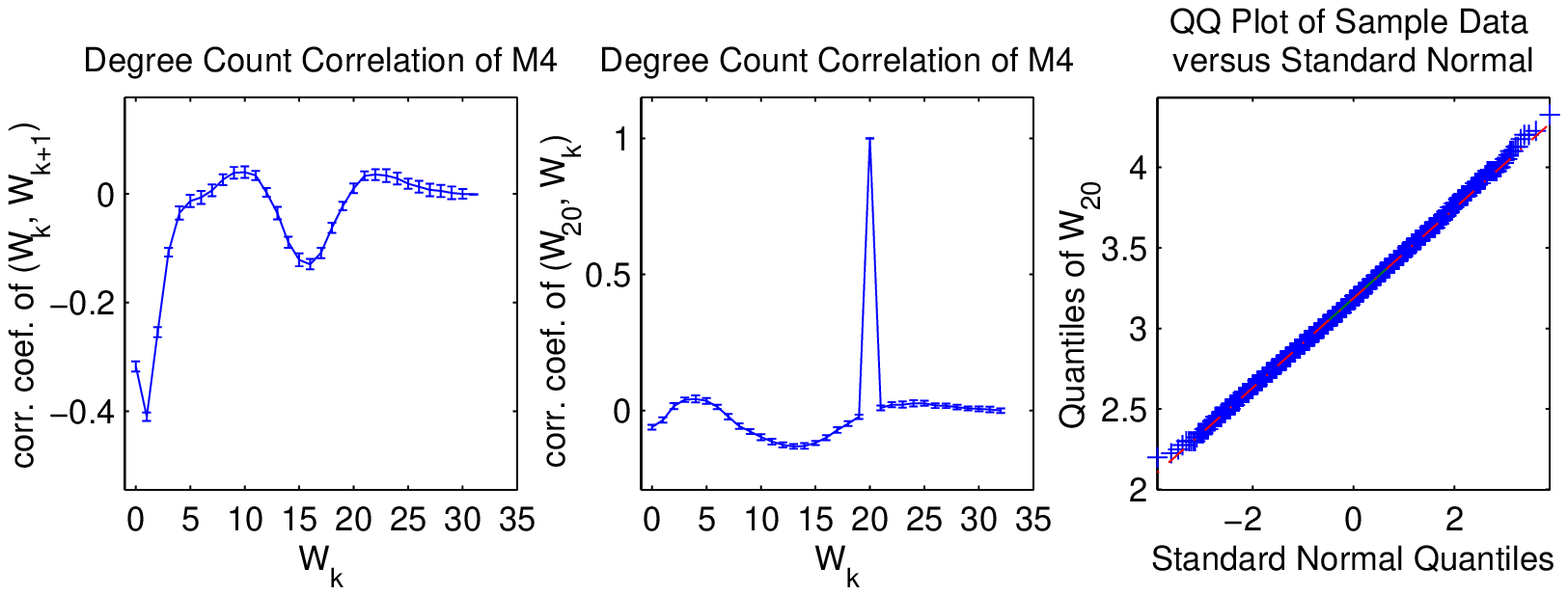}
\caption{\emph{Degree count correlation and QQ plots in Models M1 -- M4. }} \label{dccorrqq}
\end{center}
\end{figure}

Figure \ref{dccorrqq} shows the degree count correlations, firstly between degree counts for $k$ and $k+1$, and secondly
for degree counts of an asymptotically normally distributed degree count and successive degree counts; the quantile-quantile plots are given for re-assurance. We observe a strong negative correlation for degree counts which are close by, but then close to zero correlation with counts of large degrees.

\subsection{Simulations for Power-law Type Behavior}

Using Models M1 -- M4,  but now with $n=1,000$ vertices, we plot the number of vertices of degree no less than $d$ versus $d$ itself, on a log-log scale, Despite
the networks being created using independent edges, the plots seem to display a sharp linear decline, which could easily be mis-interpreted as displaying a power-law behaviour.

\begin{figure}
\begin{center}
\includegraphics[height=7cm]{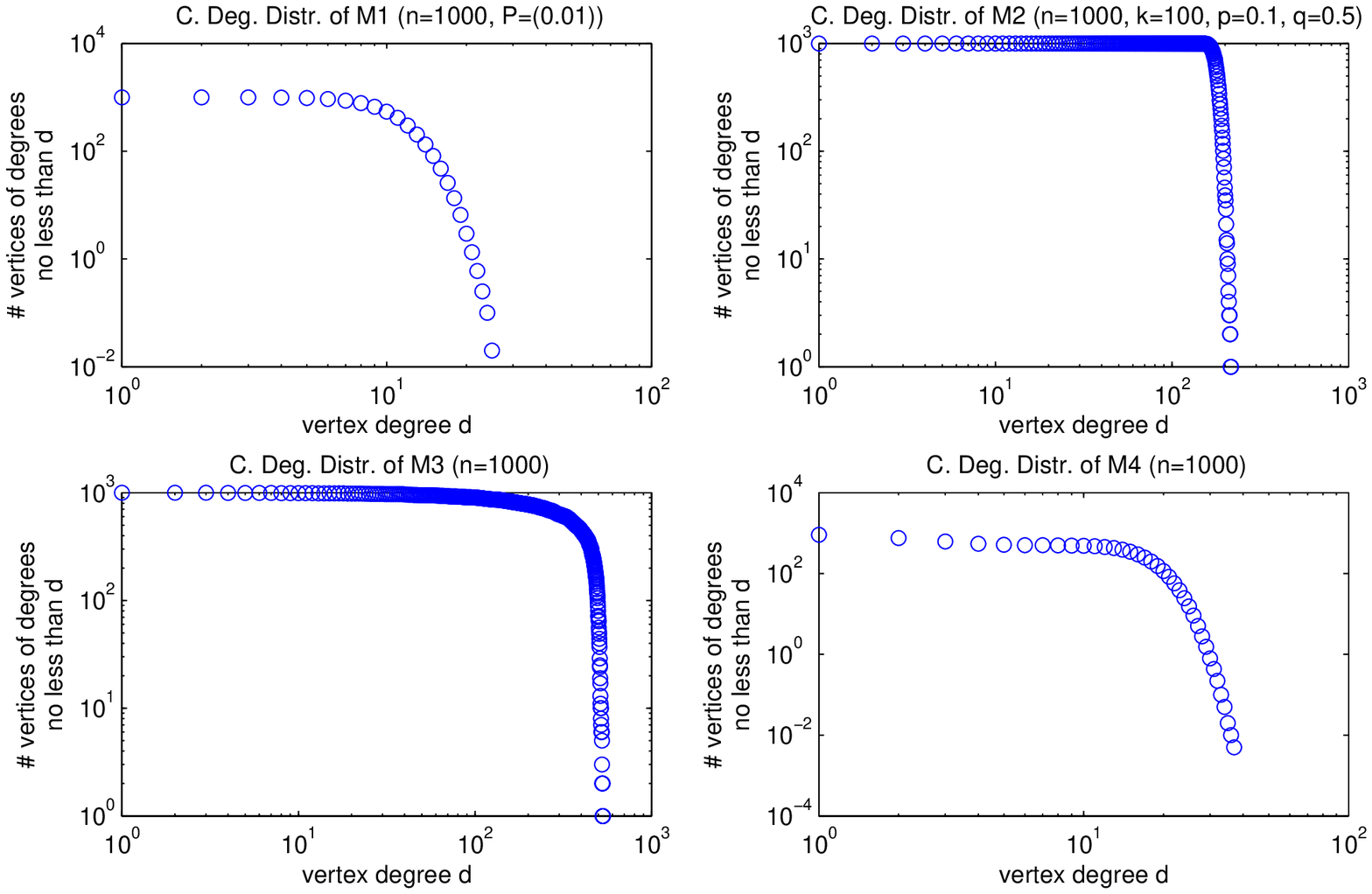}
\caption{\emph{power law; log-log scale}} \label{plawloglog}
\end{center}
\end{figure}

These simulations confirm the pseudo-power-law phenomena, and therefore raise the issue that, without rigorous analysis, simulation-based claims of detecting power-law type behaviour, or scale-free behaviour, could be in fact unreliable and misleading. The vertex degree distribution may not be a suitable visual method for distinguishing different network models. In contrast, our distributional results help assess the joint distribution of vertex degrees under a fairly general null model.

\section{Proofs} \label{proofs}

In this section, we provide proofs for the  size-biased coupling construction \ref{II}, as well as for Theorem \ref{newapprox}, Theorem \ref{mvn},   and Theorem \ref{PoissonThm}. First we prove Lemma \ref{isaprob}.

\begin{figure}
\begin{center}
\includegraphics[height=4cm]{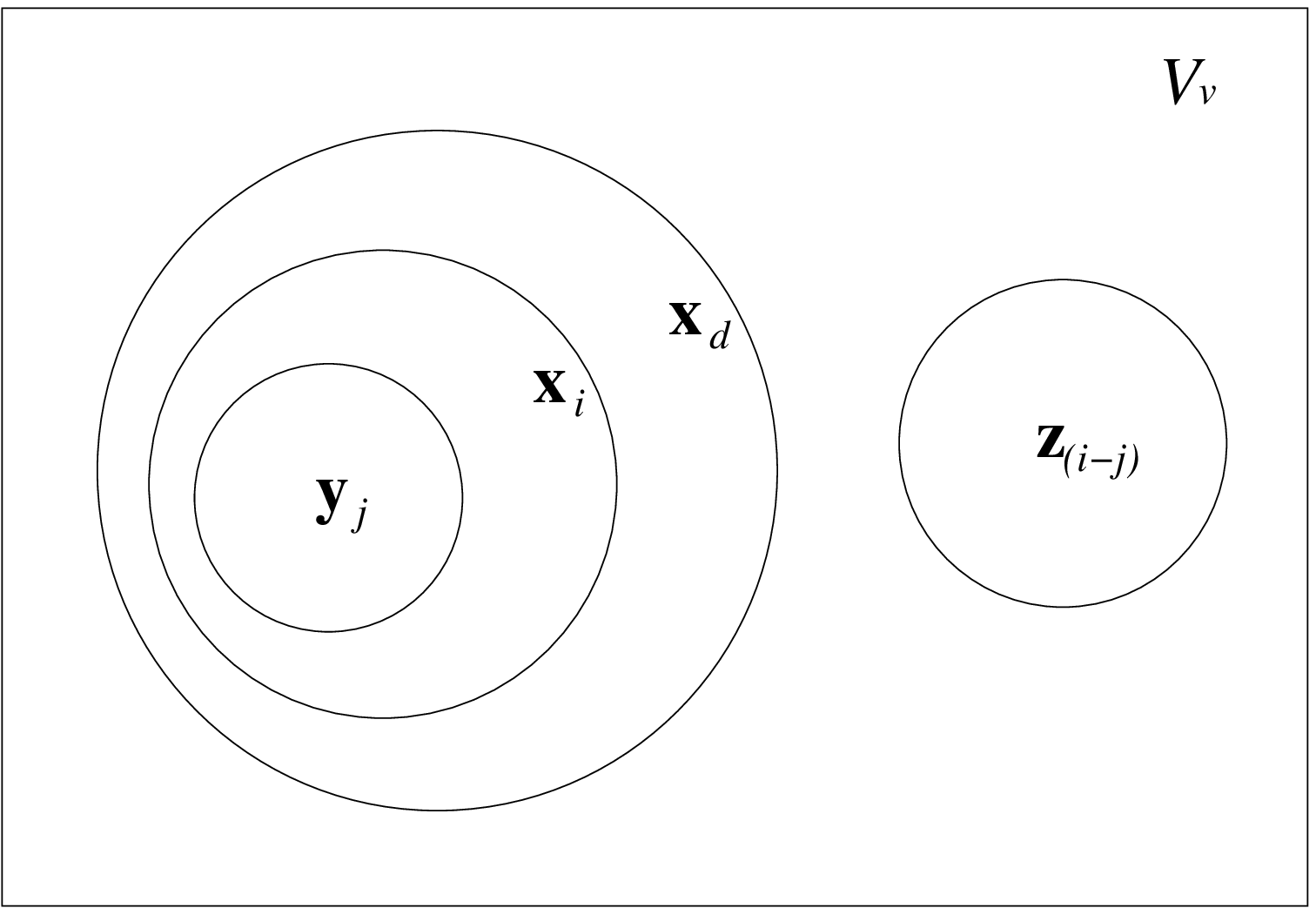}
\includegraphics[height=4cm]{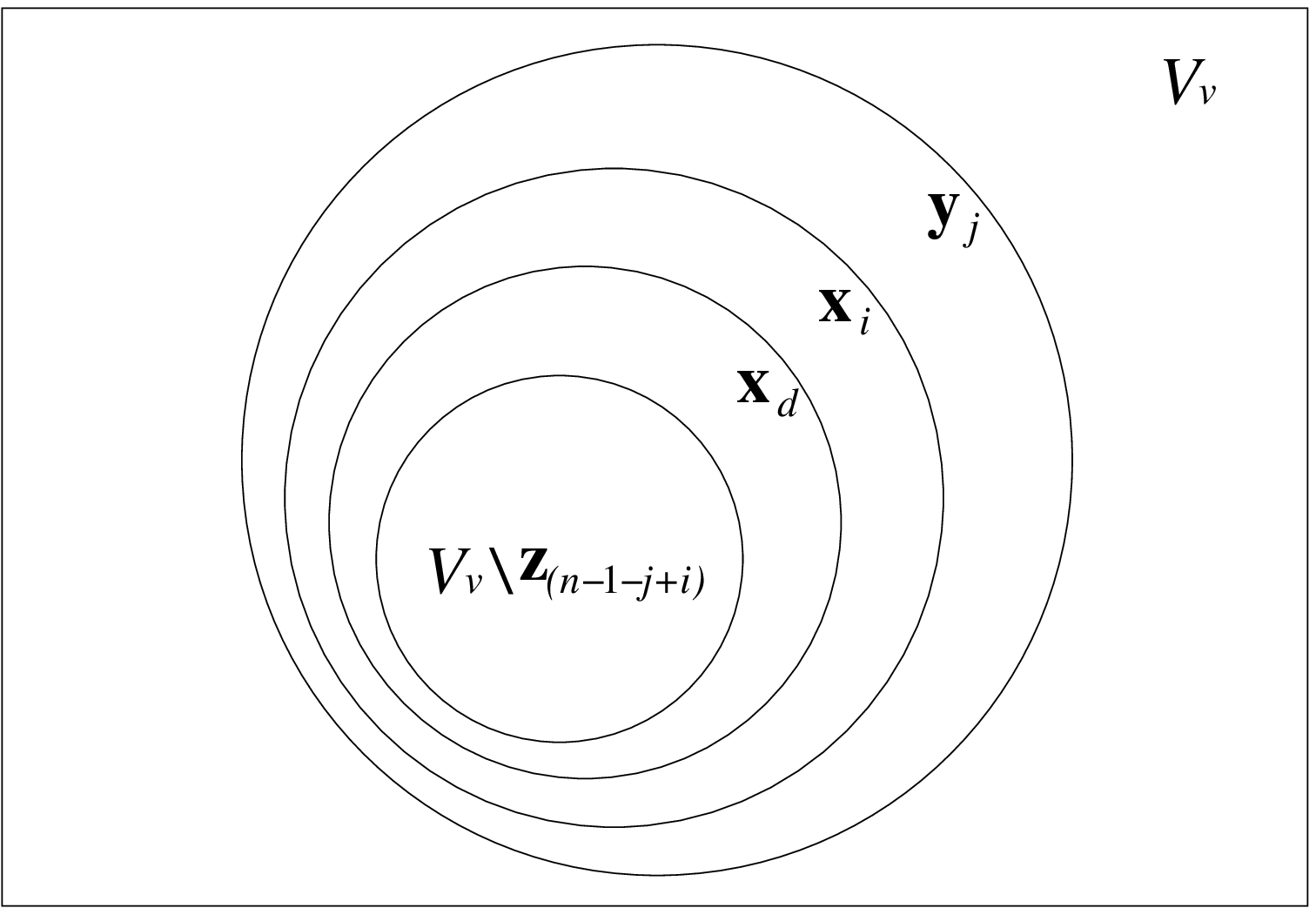}
\caption{\emph{The set diagram for the expression (\ref{delprob}) (left) and (\ref{addprob}) (right), respectively.}}\label{set12}
\end{center}
\end{figure}

\bpf{ [\textbf{Proof of  Lemma \ref{isaprob}}]} 
 We first re-write $f^{+}(\mathbf{x}_i\,|\,\mathbf{x}_d)$ as well as  $f^{-}(\mathbf{x}_i\,|\,\mathbf{x}_d)$ by writing out the set $N(v)$ in terms of those vertices which remain fixed in the construction, and those which get added or removed, respectively. Here   $\mathbf{x}_i$, $\mathbf{x}_d$, $\mathbf{y}_j$, $\mathbf{z}_{(i-j)}$ and $\mathbf{z}_{(n-1-j+i)}$ are all subsets of $V_v$.  Figure \ref{set12} illustrates the set relation in  (\ref{delprob}) and (\ref{addprob}) respectively. We have
\be
\lefteqn{f^{+}(\mathbf{x}_i\,|\,\mathbf{x}_d)} \label{delprob} \\
&=&\sum_{j=0}^{i}\,\sum_{\mathbf{y}_j:\mathbf{y}_j\subset\mathbf{x}_i}
\sum_{\stackrel{\scriptsize\mbox{$\mathbf{z}_{(i-j)}:$}}{\scriptsize\mbox{$\mathbf{z}_{(i-j)}\subset V_v\setminus\mathbf{x}_d$}}}\frac{1}{\choose{d-j}{i-j}}\,\mathbb{P}(N(v)=\mathbf{y}_j\cup\mathbf{z}_{(i-j)}\,|\,D(v)=i), \nn
\ee
and
\be
\lefteqn{f^{-}(\mathbf{x}_i\,|\,\mathbf{x}_d)} \label{addprob} \\
&=&\sum_{j=i}^{n-1} \sum_{\mathbf{y}_j:\mathbf{y}_j\supset\mathbf{x}_i}
\sum_{\stackrel{\scriptsize\mbox{$\mathbf{z}_{(n-1-j+i)}:$}}{\scriptsize\mbox{$\mathbf{z}_{(n-1-j+i)}\supset V_v\setminus\mathbf{x}_d$}}}\frac{1}{\choose{j-d}{j-i}} \mathbb{P}(N(v)=\mathbf{y}_j\cap\mathbf{z}_{(n-1-j+i)}| D(v)=i).\nn \ee

See Lemma 4.3.2 in Lin \cite{lin} for more details.

Now, from (\ref{delprob}), since, for any $\mathbf{x}_d$ and $\mathbf{y}_j$, there are {\scriptsize\mbox{$\choose{d-j}{i-j}$}} choices of $\mathbf{x}_i$, we have
\be
\lefteqn{\sum_{\mathbf{x}_i:\mathbf{x}_i\subset\mathbf{x}_d}f^{+}(\mathbf{x}_i\,|\,\mathbf{x}_d)}\nn\\
&=& \frac{1}{\mathbb{P}(D(v)=i)}\sum_{j=0}^{i}\sum_{\mathbf{y}_j:\mathbf{y}_j\subset\mathbf{x}_d}\sum_{\stackrel{\scriptsize\mbox{$\mathbf{z}_{(i-j)}:$}}{\scriptsize\mbox{$\mathbf{z}_{(i-j)}\subset V_v\setminus\mathbf{x}_d$}}}\mathbb{P}(N(v)=\mathbf{y}_j\cup\mathbf{z}_{(i-j)}).\nn
\ee

Let $\mathbf{X}_i=\{\mathbf{x}_i:\mathbf{x}_i\subset V_v\}$. Note that $\bm{\omega}_i\in\mathbf{X}_i$ if and only if $\bm{\omega}_i$ can be uniquely decomposed as $\bm{\omega}_i=\mathbf{y}_j\cup\mathbf{z}_{(i-j)}$ such that $\mathbf{y}_j\subset\mathbf{x}_d$ and $\mathbf{z}_{(i-j)}\subset V_v\setminus\mathbf{x}_d$. Thus,
\be
\sum_{\mathbf{x}_i:\mathbf{x}_i\subset\mathbf{x}_d}f^{+}(\mathbf{x}_i\,|\,\mathbf{x}_d)=
\frac{1}{\mathbb{P}(D(v)=i)}\sum_{\bm{\omega}_i\in\mathbf{X}_i}\mathbb{P}(N(v)=\bm{\omega}_i)=1,\nn
\ee
as required.
Similarly, from
 (\ref{addprob}) we find that
\be
\lefteqn{\sum_{\mathbf{x}_i:\mathbf{x}_i\supset\mathbf{x}_d}f^{-}(\mathbf{x}_i\,|\,\mathbf{x}_d)}\nn\\
&=&\frac{1}{\mathbb{P}(D(v)=i)}\sum_{j=i}^{n-1}\sum_{\mathbf{y}_j:\mathbf{y}_j\supset\mathbf{x}_d}\sum_{\stackrel{\scriptsize\mbox{$\mathbf{z}_{(n-1-j+i)}:$}}{\scriptsize\mbox{$\mathbf{z}_{(n-1-j+i)}\supset V_v\setminus\mathbf{x}_d$}}}\mathbb{P}(N(v)=\mathbf{y}_j\cap\mathbf{z}_{(n-1-j+i)}).\nn
\ee
Note that, $\bm{\omega}_i\in\mathbf{X}_i$ if and only if $\bm{\omega}_i$ can be uniquely written as $\bm{\omega}_i=\mathbf{y}_j\cap\mathbf{z}_{(n-1-j+i)}$ such that $\mathbf{y}_j\supset\mathbf{x}_d$ and $\mathbf{z}_{(n-1-j+i)}\supset V_v\setminus\mathbf{x}_d$. Thus,
\be
\sum_{\mathbf{x}_i:\mathbf{x}_i\supset\mathbf{x}_d}f^{-}(\mathbf{x}_i\,|\,\mathbf{x}_d)
=\frac{1}{\mathbb{P}(D(v)=i)}\sum_{\bm{\omega}_i\in\mathbf{X}_i}\mathbb{P}(N(v)=\bm{\omega}_i)=1,\nn
\ee
as required.
\epf

\bpf[\textbf{Proof of Construction \ref{II}}]

The goal of this proof is to show that, for $\beta= (v,i) \in\Gamma$,
\be \label{desired}
\mathcal{L}(\mathbf{X}^{\beta})=\mathcal{L}(\mathbf{X}\,|\,X_{\beta}=1),
\ee
where $\mathbf{X}=\{X_{(v,i)}:(v,i)\in\Gamma\}$ and $\mathbf{X}^{\beta}=\{X_{(v,i)}^{\beta}:(v,i)\in\Gamma\}$. Indeed we shall show that the distribution of the constructed model $\mathscr{G}^{\beta}(n,\{p_{ij}\})$ is the same as the conditional distribution of the original model $\mathscr{G}(n,\{p_{ij}\})$ given $X_{\beta}=1$, that is,
with $\mathbf{X}_{edge}=\{X_{\{a,b\}}:\{a,b\}\in E\}$ and $\mathbf{X}_{edge}^\beta=\{X_{\{a,b\}}^\beta:\{a,b\}\in E\}$,
we meed to show that for all $\mathbf{w}:=(w_{a,b},\{a,b\}\in E)\in\{0,1\}^{|E|}$,
\be
\mathbb{P}(\mathbf{X}^{(v,i)}_{edge}=\mathbf{w})=\mathbb{P}(\mathbf{X}_{edge}=\mathbf{w}\,|\,X_{(v,i)}=1).\label{RL}
\ee
The desired equation \eqref{desired}
then follows  because $\mathbf{X}$ is a function of $\mathbf{X}_{edge}$.

By definition of $X_{(v,i)}$, the right-hand side of (\ref{RL})  is zero when $\Sigma_{x\in V_v}w_{v,x}\neq i$, and by construction in that case the left-hand side of (\ref{RL})  is zero also. Assume that $\Sigma_{x\in V_v}w_{v,x} = i$, then 
the right-hand sice of (\ref{RL}) equals 
\be
\lefteqn{\mathbb{P}(\mathbf{X}_{edge}=\mathbf{w}|X_{(v,i)}=1)=\mathbb{P}(\mathbf{X}_{edge}=\mathbf{w}\,|\,D(v)=i)} \nn \\
&=&\mathbb{P}(X_{\{a,b\}}=w_{a,b}\forall a\neq v,\;b\neq v ;\;N(v)=\{x\in V_v: w_{v,x}=1\}\,|\,D(v)=i) \nn \\
&=&\mathbb{P}(X_{\{a,b\}}=w_{a,b}\;\forall \;a\neq v,\;b\neq v)\mathbb{P}(N(v)=\mathbf{x}_i\,|\,D(v)=i),\label{RHS}
\ee
where the last equality follows from  the independence of the edges, as  the condition $D(v)=i$ only affects $(X_{\{v,x\}},x\in V_v)$, but not edges which do not contain $v$.
On the other hand, the left hand side of (\ref{RL}) equals 
\be
\lefteqn{ \mathbb{P}(X_{\{a,b\}}=w_{a,b}\;\mbox{for all}\;a\neq v,\;b\neq v;\; X_{\{v,x\}}^{(v,i)}=w_{v,x}\;\mbox{for}\;x\in V_v)} \nn \\
&=&\mathbb{P}(X_{\{a,b\}}=w_{a,b}\;\mbox{for all}\;a\neq v,\;b\neq v)\mathbb{P}(X_{\{v,x\}}^{(v,i)}=w_{v,x}\;\mbox{for}\;x\in V_v),\qquad\label{LHS}
\ee
as the construction of $\mathbf{X}^{(v,i)}_{edge}$ from $\mathbf{X}_{edge}$ affects only the edges with $v$ as one of its end points, and the edges are independent.
Note that, in (\ref{LHS}),
\be
\mathbb{P}(X_{\{v,x\}}^{(v,i)}=w_{v,x}\forall x\in V_v)
 =\mathbb{P}(N^{(v,i)}(v)=\mathbf{x}_i).\nn
\ee
Hence to conclude that  (\ref{LHS}) equals (\ref{RHS}), it remains to show that
\be
\mathbb{P}(N^{(v,i)}(v)=\mathbf{x}_i)=\mathbb{P}(N(v)=\mathbf{x}_i\,|\,D(v)=i). \label{proofgoal}
\ee
Indeed,
\be
\lefteqn{\mathbb{P}(N^{(v,i)}(v)=\mathbf{x}_i)} \nn\\
&=&q_{v, i} \mathbb{P}(N(v)=\mathbf{x}_i\,|\,D(v)=i) \label{left}\\
&&+\sum_{d:d>i} q_{v, d} \sum_{\mathbf{x}_d:\mathbf{x}_d\supset\mathbf{x}_i}
\mathbb{P}(N(v)=\mathbf{x}_d | D(v)=d)\mathbb{P}(N^{(v,i)}(v)=\mathbf{x}_i | N(v)=\mathbf{x}_d)\nn\\
&&+\sum_{d:d<i} q_{v, d} \sum_{\mathbf{x}_d:\mathbf{x}_d\subset\mathbf{x}_i}
\mathbb{P}(N(v)=\mathbf{x}_d | D(v)=d)\mathbb{P}(N^{(v,i)}(v)=\mathbf{x}_i | N(v)=\mathbf{x}_d),\nn
\ee
where the three terms corresponding to the coupling construction \eqref{II}. Now, we calculate the sums over $\mathbf{x}_d:\mathbf{x}_d\supset\mathbf{x}_i$ and over $\mathbf{x}_d:\mathbf{x}_d\subset\mathbf{x}_i$ separately. In fact, for the first case $d>i$ and $\mathbf{x}_d\supset\mathbf{x}_i$, it follows from (\ref{delp}) and (\ref{delprob}) that
\be
\label{case1}
&& \sum_{\mathbf{x}_d:\mathbf{x}_d\supset\mathbf{x}_i}\mathbb{P}(N(v)=\mathbf{x}_d\,|\,D(v)=d)
\mathbb{P}(N^{(v,i)}(v)=\mathbf{x}_i\,|\,N(v)=\mathbf{x}_d) \\
&=&\sum_{j=0}^{i}\;
\sum_{\stackrel{\scriptsize\mbox{$\mathbf{x}_d:$}}{\scriptsize\mbox{$\mathbf{x}_d\supset\mathbf{x}_i$}}}\;
\sum_{\stackrel{\scriptsize\mbox{$\mathbf{y}_j:$}}{\scriptsize\mbox{$\mathbf{y}_j\subset\mathbf{x}_i$}}}\;
\sum_{\stackrel{\scriptsize\mbox{$\mathbf{z}_{(i-j)}:$}}{\scriptsize\mbox{$\mathbf{z}_{(i-j)}\subset V_v\setminus\mathbf{x}_d$}}}\frac{1}{\choose{d-j}{i-j}} \nn \\
&& \times \frac{\mathbb{P}(N(v)=\mathbf{x}_d)\mathbb{P}(N(v)=\mathbf{y}_j\cup\mathbf{z}_{(i-j)})}
{\mathbb{P}(D(v)=d)\mathbb{P}(D(v)=i)},\nn
\ee
where
\beas
\mathbb{P}(N(v)=\mathbf{x}_d)
&=&\bigg[\Big(\prod_{x\in\mathbf{x}_d}p_{vx}\Big)\Big(\prod_{x\in V_v\setminus\mathbf{x}_d} ( 1 - p_{vx}) \Big)\bigg]. 
\enas
 Since $\mathbf{y}_j\subset\mathbf{x}_i\subset\mathbf{x}_d$ and $\mathbf{z}_{(i-j)}\subset V_v\setminus\mathbf{x}_d$, we have $V_v\setminus(\mathbf{y}_j\cup\mathbf{z}_{(i-j)})\supset\mathbf{x}_d\setminus\mathbf{x}_i$ (see Figure \ref{set12} for reference). Therefore,
\be
&&{\mathbb{P}(N(v)=\mathbf{x}_d)\mathbb{P}(N(v)=\mathbf{y}_j\cup\mathbf{z}_{(i-j)})}\nn\\
&=&\bigg[\Big(\prod_{x\in\mathbf{x}_i}p_{vx}\Big)\Big(\prod_{x\in V_v\setminus\mathbf{x}_d}(1 - p_{vx}) \Big)\Big(\prod_{x\in\mathbf{x}_d\setminus\mathbf{x}_i}( 1 - p_{vx}) \Big)\bigg]\nn\\
&&\cdot\bigg[\Big(\prod_{x\in\mathbf{x}_d\setminus\mathbf{x}_i}p_{vx}\Big)
\Big(\prod_{x\in\mathbf{y}_j\cup\mathbf{z}_{(i-j)}}p_{vx}\Big)\Big(\prod_{x\in V_v\setminus[\mathbf{y}_j\cup\mathbf{z}_{(i-j)}\cup(\mathbf{x}_d\setminus\mathbf{x}_i)]}( 1 - p_{vx}) \Big)\bigg]\nn\\
&=&\mathbb{P}(N(v)=\mathbf{x}_i)\mathbb{P}\big(N(v)=\mathbf{y}_j\cup\mathbf{z}_{(i-j)}\cup(\mathbf{x}_d\setminus\mathbf{x}_i)\big).\nn
\ee
Hence, from (\ref{case1}), we have for $d>i$ that
\be
&&{\sum_{\mathbf{x}_d:\mathbf{x}_d\supset\mathbf{x}_i}\mathbb{P}(N(v)=\mathbf{x}_d\,|\,D(v)=d)
\mathbb{P}(N^{(v,i)}(v)=\mathbf{x}_i\,|\,N(v)=\mathbf{x}_d)}\nn\\
&&=\frac{\mathbb{P}(N(v)=\mathbf{x}_i)}{\mathbb{P}(D(v)=i)}\cdot\frac{1}{\mathbb{P}(D(v)=d)} \sum_{j=0}^{i} \frac{1}{\choose{d-j}{i-j}} \; \label{case1a}\\
&&
\sum_{\stackrel{\scriptsize\mbox{$\mathbf{x}_d:$}}{\scriptsize\mbox{$\mathbf{x}_d\supset\mathbf{x}_i$}}}\;
\sum_{\stackrel{\scriptsize\mbox{$\mathbf{y}_j:$}}{\scriptsize\mbox{$\mathbf{y}_j\subset\mathbf{x}_i$}}}\;
\sum_{\stackrel{\scriptsize\mbox{$\mathbf{z}_{(i-j)}:$}}{\scriptsize\mbox{$\mathbf{z}_{(i-j)}\subset V_v\setminus\mathbf{x}_d$}}}\,\mathbb{P}\big(N(v)=\mathbf{y}_j\cup\mathbf{z}_{(i-j)}\cup(\mathbf{x}_d\setminus\mathbf{x}_i)\big).\nn
\ee
Since $\mathbf{y}_j\cup\mathbf{z}_{(i-j)}$ is an $i$-set (i.e. a set with $i$ elements), $(\mathbf{x}_d\setminus\mathbf{x}_i)$ is a $(d-i)$-set, and $(\mathbf{y}_j\cup\mathbf{z}_{(i-j)})\cap(\mathbf{x}_d\setminus\mathbf{x}_i)=\emptyset$, we have that $\mathbf{y}_j\cup\mathbf{z}_{(i-j)}\cup(\mathbf{x}_d\setminus\mathbf{x}_i)$ is a $d$-set in $\mathbf{X}_d$. Conversely, for any $d$-set $\bm{\omega}_d\in\mathbf{X}_d$, since $\mathbf{x}_i$ is fixed, we can decompose $\bm{\omega}_d$ as
\be
\bm{\omega}_d=\mathbf{y}_j\cup\mathbf{z}_{(i-j)}\cup(\mathbf{x}_d\setminus\mathbf{x}_i),\label{omegad1}
\ee such that $\mathbf{y}_j\subset\mathbf{x}_i$, $\mathbf{x}_d\supset\mathbf{x}_i$ and $\mathbf{z}_{(i-j)}\subset V_v\setminus\mathbf{x}_d$. Referring to Figure \ref{set12}, it is easy to show that for any $\bm{\omega}_d\in\mathbf{X}_d$, there are {\scriptsize\mbox{$\choose{d-j}{i-j}$}} solutions $(\hat{\mathbf{y}}_j,\hat{\mathbf{x}}_d,\hat{\mathbf{z}}_{(i-j)})$ to decompose $\bm{\omega}_d$ as (\ref{omegad1}) (see Lemma 4.3.4 in Lin (2008) for more details). Thus, it follows that
\be
&&\sum_{j=0}^{i} 
\sum_{\stackrel{\scriptsize\mbox{$\mathbf{x}_d:$}}{\scriptsize\mbox{$\mathbf{x}_d\supset\mathbf{x}_i$}}}\;
\sum_{\stackrel{\scriptsize\mbox{$\mathbf{y}_j:$}}{\scriptsize\mbox{$\mathbf{y}_j\subset\mathbf{x}_i$}}}\;
\sum_{\stackrel{\scriptsize\mbox{$\mathbf{z}_{(i-j)}:$}}{\scriptsize\mbox{$\mathbf{z}_{(i-j)}\subset V_v\setminus\mathbf{x}_d$}}}\frac{1}{\choose{d-j}{i-j}}\,
\mathbb{P}\big(N(v)=\mathbf{y}_j\cup\mathbf{z}_{(i-j)}\cup(\mathbf{x}_d\setminus\mathbf{x}_i)\big)\nn\\
&&\;\;=\;\;\sum_{\bm{\omega}_d\in\mathbf{X}_d}\mathbb{P}(N(v)=\bm{\omega}_d) \;\;=\;\;\mathbb{P}(D(v)=d),\nn
\ee
and from (\ref{case1a}), we have for $d>i$ that
\be
&& { \sum_{\mathbf{x}_d:\mathbf{x}_d\supset\mathbf{x}_i} \mathbb{P}(N(v)=\mathbf{x}_d\,|\,D(v)=d)
\mathbb{P}(N^{(v,i)}(v)=\mathbf{x}_i|N(v)=\mathbf{x}_d) } \nn \\
&=&  \frac{\mathbb{P}(N(v)=\mathbf{x}_i)}{\mathbb{P}(D(v)=i)}.\quad\label{case1b}
\ee

The case $d<i$ and $\mathbf{x}_d\subset\mathbf{x}_i$ is treated similarly, giving
\be
&&{\sum_{\mathbf{x}_d:\mathbf{x}_d\subset\mathbf{x}_i}\mathbb{P}(N(v)=\mathbf{x}_d | D(v)=d)
\mathbb{P}(N^{(v,i)}(v)=\mathbf{x}_i\,|\,N(v)=\mathbf{x}_d)}\nn \\
&=& \frac{\mathbb{P}(N(v)=\mathbf{x}_i)}
{\mathbb{P}(D(v)=i)}.\quad\label{case2b}
\ee
See \cite{lin} for details. Combining (\ref{left}), (\ref{case1b}) and (\ref{case2b}),
\be
\mathbb{P}(N^{(v,i)}(v)=\mathbf{x}_i)
&=&\frac{\mathbb{P}(N(v)=\mathbf{x}_i)}{\mathbb{P}(D(v)=i)} =\mathbb{P}(N(v)=\mathbf{x}_i\,|\,D(v)=i),\nn
\ee
as required in (\ref{proofgoal}) to complete the proof.
\epf

The proof of Theorem \ref{mvn} is based on the following theorem, which is similar to Theorem 1.2 in \cite{GoldsteinRinott1996} but gives a multivariate normal approximation with respect to an alternative covariance matrix $\Sigma_0$   for which $||\Sigma_0||$ is straightforward to bound and which is close to $\Sigma$.
 We use
the notation from \cite{GoldsteinRinott1996}. For a vector $\bfb \in \R^p$ we let $\parallel \bfb \parallel = \max_{1 \le i \le p} | b_i|$. More
generally, for an array $A=(a_{i,j})$, the notation $\parallel \cdot \parallel$ is its maximal absolute value.
For an array $A( {\bf{w}})= \{ a_i({\bf{w}} ) \} $ of functions, $\parallel A \parallel= \sup_{{\bf{w}}} \max_i | a_i({\bf{w}})|$.

\bthm\label{newapprox}
{Let $\mathbf{W}=(W_i,1\leq i\leq p)$ be a random vector in $\R^p$ with nonnegative components. Let $\bm{\lambda}=(\lambda_1,\ldots,\lambda_p)=\mathbb{E}\mathbf{W}$ and assume that $Var(\mathbf{W})=\Sigma=(\sigma_{ij})$ exists. Let $\Sigma_0=(\sigma_{ij}^{0})$ be a positive definite (invertible) $p\times p$ matrix. For each $i=1,\ldots,p$ let $(\mathbf{W},\mathbf{W}^i)$ be a random vector defined on a joint probability space with $\mathbf{W}^i$ having the $\mathbf{W}$-size biased distribution in the $i^{th}$ coordinate. Let $h:\R^p\rightarrow\R\in C_b^3$, and let $Nh=\mathbb{E}h(\mathbf{Z})$ where $\mathbf{Z}$ denotes a standard normal variable in $\R^p$. Then
\beas
\lefteqn{ \left| \mathbb{E} h ( \Sigma_0^{-1/2} (\mathbf{W}-\bm{\lambda}) ) - Nh \right| }\nn\\
&\leq& \frac{p^2}{2} \|\Sigma_0^{-1/2}\|^2 \|D^2h\| \sum_{i=1}^{p} \sum_{j=1}^{p} \lambda_i \sqrt{ Var\mathbb{E} \left[ W_j^i-W_j \big| \mathbf{W} \right] }\nn\\
&& + \frac{1}{2} \frac{p^3}{3} \|\Sigma_0^{-1/2}\|^3 \|D^3h\| \sum_{i=1}^{p} \sum_{j=1}^{p} \sum_{k=1}^{p} \lambda_i \mathbb{E} \left|(W_j^i-W_j)(W_k^i-W_k)\right|\nn\\
&& + \frac{p^2}{2} \|\Sigma_0^{-1/2}\|^2 \|D^2h\| \sum_{i=1}^{p} \sum_{j=1}^{p} |\sigma_{ij}^{0}-\sigma_{ij}|.
\enas
}
\ethm

\bpf  [\textbf{Proof of Theorem \ref{newapprox}}]
The proof follows closely the lines of the proof of Theorem 1.2 in  \cite{GoldsteinRinott1996}. The only difference is that instead of their decomposition (18), we use
\bea
\lefteqn{\mathbb{E} \left\{ h( \Sigma_0^{-1/2} (\mathbf{W}-\bm{\lambda}) ) - Nh \right\}} \nn \\
&=&  -\mathbb{E} \left\{ \sum_{i,j=1}^{p} \left[\lambda_i (W_j^i-W_j)-\sigma_{ij}\right] \frac{\partial^2}{\partial w_i \partial w_j} f(\mathbf{W}) \right\} \label{Taylor1}\\
&&
-\mathbb{E} \left\{ \sum_{i,j=1}^{p}(\sigma_{ij}-\sigma_{ij}^0) \frac{\partial^2}{\partial w_i \partial w_j} f(\mathbf{W}) \right\} \label{Taylor2}\\
&&
-\mathbb{E} \left\{ \sum_{i,j,k=1}^{p} \lambda_i (W_j^i-W_j) (W_k^i-W_k) \right. \label{Taylor3}  \\
&& \left. \int_{0}^{1} (1-t) \frac{\partial^3}{\partial w_i \partial w_j \partial w_k}
f\left[ \mathbf{W} + t(\mathbf{W}^i-\mathbf{W}) \right] (W_j^i-W_j) (W_k^i-W_k) dt \right\}.  \nn
\ena
The bound for (\ref{Taylor1}) and for (\ref{Taylor3}) are as in \cite{GoldsteinRinott1996};
and for (\ref{Taylor2}), we obtain the bound
\beas
 \frac{p^2}{2} \|\Sigma_0^{-1/2}\|^2 \|D^2h\| \sum_{i=1}^{p} \sum_{j=1}^{p} |\sigma_{ij}-\sigma_{ij}^0|.
\enas
This completes the proof.
\epf

\bpf  [\textbf{Proof of Theorem \ref{mvn}}]
The proof is based on Theorem \ref{newapprox} and the size-biased coupling construction \ref{II}. Denote the randomly picked vertex by $V$. First note that, for all $i=1, \ldots, p$,
$
| W_{d_j}^{d_i} - W_{d_j} | \le |D(V) - d_i| + 1,\nn
$
because at most $|D(V) - d_i|$ edges are added or removed; and the degree of $V$ is fixed to equal $d_i$ in the size-biased distribution. Hence
\beas
\quad\mathbb{E} \left| (W_{d_j}^{d_i} - W_{d_j})(W_{d_k}^{d_i} - W_{d_k}) \right| &\le&  \mathbb{E} ( D(V) + d_i + 1)^2 \\
&\le&  2 (d_i + 1)^2 + 2 \mathbb{E} D(V)^2.\nn
\enas
Now for $d_i$ the vertex $V=v$ is chosen proportional to $ q_{v, d_i} $; hence
\beas
 \mathbb{E} D(V)^2
 &=& \frac{1}{\lambda_i} \sum_{v=1}^n  q_{v, d_i} \big(Var D(v) + (\mathbb{E} D(v))^2\big), \nn
\enas
and we use \eqref{vard} to bound the variance.

The bound on $\sqrt{ Var \mathbb{E}[W_{d_j}^{d_i} - W_{d_j} \,|\, \bfw] } $ is straightforward and follows the lines of \cite{GoldsteinRinott1996}.
First note that 
$$ Var \mathbb{E}[W_{d_j}^{d_i} - W_{d_j} \,|\, \bfw] \le  Var \mathbb{E}[W_{d_j}^{d_i} - W_{d_j} \,|\, \mathscr{G}(n, \{ p_{ij} \}) ].$$ 
With the notation \eqref{delp} and \eqref{addp}, we abbreviate
\beas
&& \alpha_{i,j}(u,v) \nn \\
&=& \alpha(u,v) \\
&=& \frac{q_{v,d_i}}{\lambda_i}  \sum_{d > d_i}  \sum_{ \mathbf{x}_d: u \in \mathbf{x}_d} \mathbbm{1}(N(v) = \mathbf{x}_d)
 \left\{ \mathbbm{1}(D(u) = d_j + 1) - \mathbbm{1}(D(u) = d_j) \right\} \nn\\
&& \times \sum_{ \mathbf{x}_{d_i}\subset\mathbf{x}_d; u \not\in \mathbf{x}_{d_i}} f^{+}(\mathbf{x}_{d_i}\,|\,\mathbf{x}_d) \nn
\enas
and
\bea
&& \beta_{i,j}(u,v) \nn \\
&=& \beta(u,v) \label{beta} \\
&=& \frac{q_{v,d_i}}{\lambda_i} \sum_{d < d_i}  \sum_{ \mathbf{x}_d: u \not\in \mathbf{x}_d} \mathbbm{1}(N(v) = \mathbf{x}_d)
 \left\{ \mathbbm{1}(D(u) = d_j - 1) - \mathbbm{1}(D(u) = d_j) \right\} \nn\\
&& \times \sum_{ \mathbf{x}_{d_i}\supset\mathbf{x}_d; u \in \mathbf{x}_{d_i}} f^{-}(\mathbf{x}_{d_i}\,|\,\mathbf{x}_d). \nonumber
\ena
We note that
\bea \label{alphabound}
\vert \alpha(u,v) \vert \leq  \frac{q_{v,d_i}}{\lambda_i} \quad\mbox{ and }\quad \vert \beta(u,v) \vert \leq  \frac{q_{v,d_i}}{\lambda_i} .
\ena
Then
\beas
\lefteqn{\mathbb{E} \big[W_{d_j}^{d_i} - W_{d_j} \,\big|\, \mathscr{G}(n, \{ p_{ij} \}) \big] }\nn\\
 &=&  \sum_v \sum_{u \ne v} \mathbbm{1}(u \sim v) \alpha(u,v)  +  \sum_v \sum_{u \ne v} \mathbbm{1}(u \not\sim v) \beta (u,v)\nn\\
&&+ \frac{1}{\lambda_i} \sum_v  q_{v, d_i}\mathbbm{1}(D(v)\neq d_i)\mathbbm{1}(i=j) - \frac{1}{\lambda_i} \sum_v  q_{v, d_i} \mathbbm{1}(D(v)=d_j)
\mathbbm{1}(i \ne j).\nn
\enas
This gives that
\bea
&&  Var\mathbb{E}\big[ W_{d_j}^{d_i} - W_{d_j} \,\big|\, \bfw \big] \nn \\
&\leq& 4 \left\{ Var \left(  \sum_v   \sum_{u \ne v} \mathbbm{1}(u \sim v) \alpha(u,v) \right) \right.  \label{term1}\\
&&\quad\left. +  Var \left(  \sum_v   \sum_{u \ne v} \mathbbm{1}(u \not\sim v) \beta(u,v) \right)
 \right. \label{term2} \\
&&\quad\left. +  Var \left( \frac{1}{\lambda_i} \sum_v  q_{v, d_i}\mathbbm{1}(D(v)\neq d_i)\mathbbm{1}(i=j) \right)
 \right. \label{term3} \\
&&\quad\left. +  Var \left(  \frac{1}{\lambda_i} \sum_v  q_{v, d_i} \mathbbm{1}(D(v)=d_j)
\mathbbm{1}(i \ne j)  \right) \right\} \qquad \label{term4}.
\ena
Firstly, for \eqref{term1},
with \eqref{alphabound}, and $u' \ne u, v$ as well as $ v' \ne u, v$, 
\beas
&&Var \big(  \mathbbm{1}(u \sim v) \alpha(u,v)  \big) \le \mathbb{E}  \big(  \mathbbm{1}(u \sim v) \alpha(u,v) \big)^2
\le  p_{u,v}\frac{q_{v,d_i}^2}{\lambda_i^2}; \nn \\
&&\big| Cov \big(  \mathbbm{1}(u \sim v) \alpha(u,v) ;  \mathbbm{1}(u' \sim v) \alpha(u',v) \big) \big|
\le 2 p_{u,v}  p_{u',v}\frac{q_{v,d_i}^2}{\lambda_i^2} ; \nn \\
&&\big| Cov \big(  \mathbbm{1}(u \sim v) \alpha(u,v) ;  \mathbbm{1}(u \sim v') \alpha(u,v') \big)\big|
\le 2  p_{u,v}  p_{u,v'}\frac{q_{v,d_i}q_{v',d_i}}{\lambda_i^2}.\nn
\enas
For $u'\ne u, v$ and $v'\ne u, u', v$, let $C = C (u,u',v,v') = \mathbbm{1}( u \not\sim u', u \not\sim v', v \not\sim v', u' \not\sim v)$.
Then
\bea \label{csetbound}
\mathbb{P}(C(u,u',v,v') =1 )
&\ge& 1 - (p_{u,u'} +  p_{u,v'} + p_{v,v'} + p_{u',v}).
\ena
Put $\lambda(u,v) = \mathbb{E} \mathbbm{1}(u \sim v) \alpha(u,v) $; then for $u,u',v,v'$ mutually different,
\beas
\lefteqn{ Cov\big( \mathbbm{1}(u \sim v) \alpha(u,v) , \mathbbm{1}(u' \sim v')\alpha(u',v') \big)} \nn\\
&=&
\mathbb{E} \big[ (\mathbbm{1}(u \sim v) \alpha(u,v) - \lambda(u,v) ) (\mathbbm{1}(u' \sim v') \alpha(u',v') - \lambda(u',v') ) \big\vert C=1 \big]\nn\\
&& \mathbb{P}(C=1) \nn\\
&& +
\mathbb{E} \big[(\mathbbm{1}(u \sim v) \alpha(u,v) - \lambda(u,v)) (\mathbbm{1}(u' \sim v') \alpha(u',v') - \lambda(u',v') ) \big\vert C=0 \big]\nn\\
&& \mathbb{P}(C=0).\nn
\enas
As $ \mathbbm{1}(u \sim v) | \alpha(u,v) | \le \mathbbm{1}(u \sim v)  \frac{q_{v,d_i}}{\lambda_i}$ by \eqref{alphabound} and as the  edge indicators are independent, we can bound
\beas
&&  \Big| \mathbb{E} \big[ ( \mathbbm{1}(u \sim v) \alpha(u,v)  - \lambda(u,v) ) ( \mathbbm{1}(u' \sim v') \alpha(u',v')
- \lambda(u',v') ) \big| C =0 \big] \Big|\nn\\
&& \mathbb{P}(C=0) \nn\\
&& \quad\leq 4 \frac{q_{v,d_i}q_{v',d_i} }{\lambda_i^2}  p_{u,v} p_{u', v'} (p_{u,u'} +  p_{u,v'} + p_{v,v'} + p_{u',v}) .\nn
\enas
As $\alpha(u,v)$ is a random variable which depends only on
$\{ \mathbbm{1}(u \sim v), \mathbbm{1}(w \sim u), \mathbbm{1}(w \sim v), w \ne u,v\}$, it follows that conditional on $C=1$, $\alpha(u,v) $ and $\alpha(u',v')$ are independent.
Moreover,
\beas
\lefteqn{ \mathbb{E} \big[ \mathbbm{1}(u \sim v) \alpha(u,v)  - \lambda(u,v) \big\vert C =1\big] }\nn\\
&=& \frac{q_{v,d_i}}{\lambda_i} \sum_{d > d_i}  \sum_{ \mathbf{x}_d: u \in \mathbf{x}_d} \sum_{ \mathbf{x}_{d_i}\subset\mathbf{x}_d; u \not\in \mathbf{x}_{d_i}} f^{+}(\mathbf{x}_{d_i}\,|\,\mathbf{x}_d) \nn \\
&& \times \Big\{ \mathbb{E} \big[ \mathbbm{1}(N(v) = \mathbf{x}_d) \left\{\mathbbm{1}(D(u) = d_j + 1) - \mathbbm{1}(D(u) = d_j)\right\}
\big\vert C=1\big]\nn\\
&& \quad - \mathbb{E} \big[ \mathbbm{1}(N(v) = \mathbf{x}_d) \left\{ \mathbbm{1}(D(u) = d_j + 1) - \mathbbm{1}(D(u) = d_j) \right\} \big] \Big\}. \nn
\enas
Re-grouping the terms and conditioning on $u\sim v$ give that
 \beas
\lefteqn{ \mathbb{E}\big[ \mathbbm{1}(N(v) = \mathbf{x}_d) \mathbbm{1}(u \sim v) \mathbbm{1}(D(u) = d_j + 1)  \big\vert C =1\big] }  \nn\\
 && - \mathbb{E} \big[ \mathbbm{1}(N(v) = \mathbf{x}_d)  \mathbbm{1}(u \sim v) \mathbbm{1}(D(u) = d_j + 1)\big]  \nonumber \\
 &=& p_{u,v} \Big\{ \big[ \mathbb{P} (N^{(u,u',v')} (v) = \mathbf{x}_d  \setminus \{ u\} )  -  \mathbb{P}(N^{(u)}(v) = \mathbf{x}_d \setminus \{ u\} ) \big] \nn\\
&& \mathbb{P}( D^{(u',v,v')}(u) = d_j )  \nonumber \\
 &&  - \mathbb{P} (N^{(u)}(v) = \mathbf{x}_d \setminus \{ u\} ) \big[ \mathbb{P}( D^{(v)} (u) = d_j ) - \mathbb{P}( D^{(u',v,v')}(u) = d_j )\big]\Big\}. \nonumber
\enas
Now, conditioning on whether or not $u' \sim v$ and $v' \sim v$ gives
that
\beas
\lefteqn{\Big| \mathbb{E}\mathbbm{1}(N^{(u,u',v')} (v) = \mathbf{x}_d  \setminus \{ u\} ) - \mathbb{E}\mathbbm{1}(N^{(u)}(v) = \mathbf{x}_d \setminus \{ u\} ) \Big|} \\
 &\leq & (p_{u',v} + p_{v',v}) \left\{ \mathbb{P}(N^{(u,u',v')}(v) = \mathbf{x}_d  \setminus \{ u\} )\right. \\
&&\left.
+  \mathbb{P}(N^{(u, u', v')}(v) = \mathbf{x}_d \setminus \{ u, u'\}) \right. \\
&&\left. + \mathbb{P}(N^{(u, u', v')}(v) = \mathbf{x}_d \setminus \{ u, v'\}) + \mathbb{P}(N^{(u, u', v')}(v) = \mathbf{x}_d \setminus \{ u, u', v'\})  \right\}.
\enas
Similarly,
$
\Big\vert \mathbb{P}( D^{(v)} (u) = d_j ) - \mathbb{P}( D^{(u',v,v')}(u) = d_j ) \Big\vert \le  p_{u,u'} + p_{u,v'}. \nn
$
Hence
\beas
&& \Big\vert \mathbb{E}\big[ \mathbbm{1}(N(v) = \mathbf{x}_d)\mathbbm{1}(u \sim v)\mathbbm{1}(D(u) = d_j + 1)  \big\vert C =1\big] \nn\\
&&  - \mathbb{E} \big[ \mathbbm{1}(N(v) = \mathbf{x}_d)  \mathbbm{1}(u \sim v) \mathbbm{1}(D(u) = d_j + 1)\big] \Big\vert \nonumber \\
&\leq&    p_{u,v} \left\{  (p_{u',v}+p_{v',v}) \mathbb{P}( N^{(u, u', v')}(v)  \in A ) \right.\\
&&\left. +  (p_{u,u'} + p_{u,v'}) \mathbb{P}(N^{(u)}(v) = \mathbf{x}_d  \setminus \{ u\} ) \right\} ,
\enas
where  $A=\{ \mathbf{x}_d  \setminus \{ u\}, \mathbf{x}_d \setminus \{ u, u'\}, \mathbf{x}_d \setminus \{ u, v'\},  \mathbf{x}_d \setminus \{ u, u', v'\}\}$.
%
Combining these bounds,
\beas
&& \Big\vert \mathbb{E} \big[  \big(  \mathbbm{1}(u \sim v) \alpha(u,v)  - \lambda(u,v) \big)
\big( \mathbbm{1}(u' \sim v') \alpha(u',v') - \lambda(u',v') \big)  \big\vert C=1\big] \Big\vert \nonumber\\
&\leq & \;\; 64 \, p_{u,v} p_{u',v'} ( p_{u,u'}+p_{u,v'}+p_{u',v} + p_{v',v})^2 \,\frac{ q_{v,d_i} q_{v',d_i} }{\lambda_i^2} .  \nonumber
\enas

\medskip
A similar argument holds for \eqref{term2}, involving $ \mathbbm{1}(u \not\sim v) \beta (u,v)$; recall \eqref{beta}.
Firstly,
\beas
\lefteqn{ Var \big( \mathbbm{1}(u \not\sim v) \beta (u,v)  \big)}\nn\\
&\leq& \frac{q_{v,d_i}^2}{\lambda_i^2} \mathbb{E} \bigg( \sum_{d < d_i}  \sum_{ \mathbf{x}_d: u \not\in \mathbf{x}_d} \mathbbm{1}(N(v) = \mathbf{x}_d)
 \left\{ \mathbbm{1}(D(u) = d_j - 1) - \mathbbm{1}(D(u) = d_j) \right\} \nn\\
 &&\qquad\qquad \times \;  \mathbbm{1}(u  \not\in \mathbf{x}_d) \sum_{ \mathbf{x}_{d_i}\supset \mathbf{x}_d; u \in \mathbf{x}_{d_i}} f^{-}(\mathbf{x}_{d_i}\,|\,\mathbf{x}_d)
  \bigg)^2 .\nn
\enas
We bound the probability that vertex $u$ is picked to be added to the neighbours of $v$, if $N(v) = \mathbf{x}_d$,
\bea
\mathbb{P}(u \mbox{ picked}\,|\, N(v) = \mathbf{x}_d)
&=&  \mathbbm{1}(u  \not\in \mathbf{x}_d) \sum_{ \mathbf{x}_{d_i}\supset \mathbf{x}_d; u \in \mathbf{x}_{d_i}} f^{-}(\mathbf{x}_{d_i}\,|\,\mathbf{x}_d) \\
&\leq&  p_{u,v} \frac{q_{v, d_i-1}^{(u)}}{q_{v,d_i}} \label{upickedbound}.
\ena
With \eqref{upickedbound} we obtain that
\bea
{  \sum_v   \sum_{u \ne v} Var \big( \mathbbm{1}(u \not\sim v) \beta (u,v)  \big)}
 &\leq&  \sum_v   \sum_{u \ne v}  p_{u,v}^2 \left( \frac{q_{v,d_i-1}^{(u)}}{\lambda_i} \right)^2 \label{beta1bound} .
\ena

\bigskip
Similarly as above,
 we obtain
\bea \label{beta2bound}
\lefteqn{ \sum_v \sum_{u\neq v} \sum_{u'\neq u,v} Cov\big( \mathbbm{1}(u \not\sim v) \beta (u,v) , \mathbbm{1}(u' \not\sim v) \beta (u',v) \big) } \nn\\
&\leq&
2\sum_v \sum_{u\neq v} \sum_{u'\neq u,v} p_{u,v} p_{u',v} \frac{q_{v, d_i-1}^{(u)} q_{v, d_i-1}^{(u')} }{\lambda_i^2}
\ena
and
\bea \label{beta3bound}
\lefteqn{ \sum_v \sum_{u\neq v} \sum_{v'\neq u,v} Cov\big( \mathbbm{1}(u \not\sim v) \beta (u,v) , \mathbbm{1}(u \not\sim v') \beta (u,v')  \big) } \nn \\
&\leq& 2\sum_v \sum_{u\neq v} \sum_{v'\neq u,v}  p_{u,v} p_{u,v'} \frac{q_{v, d_i-1}^{(u)} q_{v', d_i-1}^{(u)} }{\lambda_i^2}  .
\ena

\bigskip
Now assume that $u, v, u'$ and $v'$ are all distinct. We refine the definition of $\beta$; for $t=1,2$, define $\beta^{(t)}$ by
\beas
{  \beta^{(1)} (u,v)}
&=& \frac{q_{v,d_i}}{\lambda_i} \sum_{d < d_i}  \sum_{ \mathbf{x}_d: u \not\in \mathbf{x}_d} \mathbbm{1}(N(v) = \mathbf{x}_d; D(u) = d_j - 1) \\
&& 
\sum_{ \mathbf{x}_{d_i}\supset\mathbf{x}_d; u \in \mathbf{x}_{d_i}} f^{-}(\mathbf{x}_{d_i}|\mathbf{x}_d)  \nn\\
  \beta^{(2)} (u,v)
&=& \frac{q_{v,d_i}}{\lambda_i} \sum_{d < d_i}  \sum_{ \mathbf{x}_d: u \not\in \mathbf{x}_d} \mathbbm{1}(N(v) = \mathbf{x}_d; D(u) = d_j) \nn \\
&& 
\sum_{ \mathbf{x}_{d_i}\supset\mathbf{x}_d; u \in \mathbf{x}_{d_i}} f^{-}(\mathbf{x}_{d_i}|\mathbf{x}_d).\nn
\enas
Now,
\beas
\lefteqn{ \mathbb{E} \mathbbm{1}(u \not\sim v) \beta^{(1)} (u,v) } \nn\\
&=& \frac{q_{v,d_i}}{\lambda_i} \sum_{d < d_i}  \sum_{ \mathbf{x}_d: u \not\in \mathbf{x}_d}  \Big\{ p_{u,v} \mathbb{P}(D^{(v)} (u) = d_j - 2)
\mathbb{P}( N^{(u)}(v) = \mathbf{x}_d\setminus \{ u\} ) \nn\\
&&  + (1 - p_{u,v})  \mathbb{P}(D^{(v)}(u) = d_j - 1) \mathbb{P}(
N^{(u)}(v) = \mathbf{x}_d) \Big\} \\
&& \mathbb{P}( u \mbox{ picked} \, | \, N(v) =  \mathbf{x}_d)\nn
\enas
and
\beas
\lefteqn{ \mathbb{E} \mathbbm{1}(u \not\sim v) \beta^{(1)} (u,v)  \mathbbm{1}(u' \not\sim v') \beta^{(1)} (u',v')  }\nn\\
&=&  \frac{q_{v,d_i}q_{v',d_i}}{\lambda_i^2} \sum_{d < d_i}  \sum_{d' < d_i} \sum_{ \mathbf{x}_d: u \not\in \mathbf{x}_d}
\sum_{ \mathbf{y}_{d'}: u \not\in\mathbf{y}_{d'}} \mathbb{P} ( u \mbox{ picked} \,|\,   N(v) = \mathbf{x}_d )   \nn\\[1ex]
&&  \times\; \mathbb{P} ( u' \mbox{ picked} \,|\,    N(v') = \mathbf{y}_{d'} ) \nn \\
&& \times \mathbb{P}(   N(v) = \mathbf{x}_d,   N(v') = \mathbf{y}_{d'} , D(u) = d_j - 1, D(u') = d_j - 1).\nn
\enas
Moreover, when conditioning on  $C=1$, the independence of the edges gives
\beas
\lefteqn{  \mathbb{P}\big(   N(v) = \mathbf{x}_d,   N(v') = \mathbf{y}_{d'} , D(u) = d_j - 1, D(u') = d_j - 1 \big| C=1\big) } \nn\\
&=&  \mathbb{P}(   N^{(u,u',v')} (v) = \mathbf{x}_d)\mathbb{P}(   N^{(u,v,u')} (v') = \mathbf{y}_{d'} ) \nn\\
&& \mathbb{P}( D^{(v,u',v')} (u) = d_j - 1) \mathbb{P}(D^{(u,v,v')} (u') = d_j - 1).\nn
\enas
Conditioning on whether or not $u \sim v$ and $u' \sim v'$ yields
 \beas
\lefteqn{ Cov\big( \mathbbm{1}(u \not\sim v) \beta^{(1)} (u,v) , \mathbbm{1}(u' \not\sim v') \beta^{(1)} (u',v')  \big) } \nn\\
&=& \frac{q_{v,d_i}q_{v',d_i}}{\lambda_i^2} \sum_{d < d_i}  \sum_{d' < d_i} \sum_{ \mathbf{x}_d: u \not\in \mathbf{x}_d}
\sum_{ \mathbf{y}_{d'}: u \not\in\mathbf{y}_{d'}}   \mathbb{P} ( u \mbox{ picked} \,|\,   N(v) = \mathbf{x}_d ) \nn \\[1ex]
&& \times
\mathbb{P} ( u' \mbox{ picked} \,|\,   N(v') = \mathbf{y}_{d'} )\nn \\
&& \Big\{   \mathbb{P}(   N^{(u,u',v')} (v) = \mathbf{x}_d)\mathbb{P}(   N^{(u,v,u')} (v') = \mathbf{y}_{d'} )\mathbb{P}( D^{(v,u',v')} (u) = d_j - 1) \nn \\
&& \times \mathbb{P}(D^{(u,v,v')} (u') = d_j - 1) -  \mathbb{P} (D^{(v)}(u) = d_j - 1)\mathbb{P}(  N^{(u)}(v) = \mathbf{x}_d  )  \nn\\
&&\times \mathbb{P} (D^{(v')}(u') = d_j - 1)\mathbb{P}( N^{(u')} (v') = \mathbf{y}_{d'} ) \Big\} + R_1 + R_2,\nn
\enas
where
from \eqref{upickedbound} we immediately get
\beas
| R_1|  &\leq&  \max\{ p_{u,v}, p_{u',v'} \} p_{u,v}p_{u',v'} \frac{q_{v, d_i-1}^{(u)} q_{v', d_i-1}^{(u')}}{\lambda_i^2}.\nn
\enas
and, with \eqref{csetbound},
\beas
| R_2 |
&=& (p_{u,u'} + p_{u,v'} + p_{v,v'} + p_{u',v} ) p_{u,v}p_{u',v'} \frac{q_{v, d_i-1}^{(u)} q_{v', d_i-1}^{(u')}}{\lambda_i^2} .\nn
\enas
Again conditioning on the presence of edges,
we obtain
\beas
\lefteqn{ \left| \mathbb{P}(   N^{(u,u',v')} (v) = \mathbf{x}_d) -  \mathbb{P}(   N^{(u)} (v) = \mathbf{x}_d) \right| }\nn\\
&\leq&
 \mathbbm{1}( \{ u', v'\} \cap \mathbf{x}_d \ne \emptyset)  \mathbb{P}(   N^{(u)} (v) = \mathbf{x}_d)\nn\\
&& +  p_{v',v} \mathbb{P}(   N^{(u,u',v')} (v) = \mathbf{x}_d; u' \not\sim v\,|\,  v' \sim v  ) \nn\\
&& + p_{u',v} \left\{ \mathbb{P}(   N^{(u,u',v')} (v) = \mathbf{x}_d; v' \not\sim v  \,|\,  u' \sim v ) \right.  \nn\\
&&\left. +  p_{v',v} \mathbb{P}(   N^{(u,u',v')} (v) = \mathbf{x}_d \,|\,  u' \sim v, v' \sim v  ) \right\} .\nn
\enas
We also have that
\beas
 \left| \mathbb{P} ( D^{(v)} (u) = d_j - 1; u \not\sim u', u \not\sim v'  )  - \mathbb{P} ( D^{(v)} (u) = d_j - 1 ) \right|
 &\leq & p_{u,u'} + p_{u,v'}.\nn
\enas
With similar bounds for the other two terms we note that the sums over $\mathbf{x}_d$ and $ \mathbf{y}_{d'}$ still include a term of the form $ \mathbb{P}(   N^{(u,u',v')} (v) = \mathbf{x}_d; u' \not\sim v\,|\,  v' \sim v  ) $ or similar, so that these sums can be bounded by 1. Moreover, 
$$
 \sum_{ \mathbf{x}_d: u\not\in\mathbf{x}_d} \mathbbm{1}( \{ u', v'\} \cap \mathbf{x}_d \ne \emptyset)  \mathbb{P}(   N^{(u)} (v) = \mathbf{x}_d)
 \leq  p_{u',v} + p_{v',v}.  $$
We conclude that
\beas
\lefteqn{ \left|  Cov\big( \mathbbm{1}(u \not\sim v) \beta^{(1)} (u,v) , \mathbbm{1}(u' \not\sim v') \beta^{(1)} (u',v')  \big)  \right| } \nn\\
&\leq&   (p_{u,v}+p_{u',v'}+p_{u,u'} + p_{u,v'} + p_{v,v'} + p_{u',v} )  p_{u,v}p_{u',v'} \frac{q_{v, d_i-1}^{(u)} q_{v', d_i-1}^{(u')}}{\lambda_i^2} \\&& ( d_i^2 + 2)  .\nn
\enas
In  the same way we can bound
$
\left|  Cov\big( \mathbbm{1}(u \not\sim v) \beta^{(1)} (u,v) , \mathbbm{1}(u' \not\sim v') \beta^{(2)} (u',v')  \big)  \right| $.
Thus we obtain that
\bea \label{beta4bound}
\lefteqn{ \sum_v   \sum_{u \ne v} \sum_{u' \ne u,v} \sum_{v' \ne u, u', v} Cov\big( \mathbbm{1}(u \not\sim v) \beta (u,v) , \mathbbm{1}(u' \not\sim v') \beta (u',v') \big)} \nonumber \\
&\leq& 4  ( d_i^2 + 2) \sum_v   \sum_{u \ne v}  \sum_{u' \ne u,v} \sum_{v' \ne u, u', v} \frac{q_{v, d_i-1}^{(u)} q_{v', d_i-1}^{(u')}}{\lambda_i^2} \nonumber \\
&&\times\; (p_{u,v}+p_{u',v'}+p_{u,u'} + p_{u,v'} + p_{v,v'} + p_{u',v} ) p_{u,v}p_{u',v'}  .
\ena

\medskip
For \eqref{term3} use the bound
\beas
\lefteqn{ Var \left( \frac{1}{\lambda_i} \sum_v  q_{v, d_i}\mathbbm{1}(D(v)\neq d_i) \mathbbm{1}(i=j)\right) }\nn\\
&\le& \frac{1}{\lambda_i^2 } \sum_v  q_{v, d_i}^2  + \frac{1}{\lambda_i^2 } \sum_v \sum_{u \ne v}  q_{v, d_i} q_{u, d_i} Cov\big( \mathbbm{1}(D(v) \neq d_i), \mathbbm{1}(D(u) \neq d_i)\big)\nn .
\enas
Using that
\bea \label{qdiff}
|  q^{(u)}_{v, d_i } - q_{v, d_i} |
&\leq& \mathbb{P} (  D^{(u)}(v)  \neq D(v)) =  p_{u,v} ,
\ena
we obtain for \eqref{term3} that, when $i=j$,
 \bea
{ Var \left( \frac{1}{\lambda_i} \sum_v  q_{v, d_i}\mathbbm{1}(D(v)\neq d_i)  \right) }
&\le& \frac{1}{\lambda_i^2 } \sum_v  q_{v, d_i}^2 + \frac{3}{\lambda_i^2 } \sum_v \sum_{u \ne v}  q_{v, d_i} q_{u, d_i} p_{u,v} \nn\\
&\leq& \frac{3}{\lambda_i^2 } \sum_v \sum_u q_{v, d_i} q_{u, d_i} p_{u,v},\label{term3bound}
\ena
where we use the convention $p_{v,v}=1$. Similarly, for \eqref{term4}, when $i \ne j$,
\bea
 Var \left( \frac{1}{\lambda_i} \sum_v  q_{v, d_i} \mathbbm{1}(D(v)=d_j) \right)
&\le &   \frac{3}{\lambda_i^2 } \sum_v \sum_u q_{v, d_i} q_{u, d_i} p_{u,v} \label{term4bound}.
\ena
Combining the bounds for \eqref{term1} with  \eqref{term3bound},  \eqref{beta1bound}, \eqref{beta2bound}, \eqref{beta3bound}, and \eqref{beta4bound}, \eqref{qdiff} and \eqref{term4bound}, and using crude bounds such as $q_{v, d_i} \le 1$, we obtain that
\beas
Var\mathbb{E}\big[ W_{d_j}^{d_i} - W_{d_j} \,\big|\, \bfw \big] &\le& B_i^2/\lambda_i^2
\enas
with $B_i$ given in the statement of Theorem \ref{mvn}.

\medskip
The next step is  bounding $\parallel  \Sigma_0^{-1/2} \parallel $. Define
$
b_v(i)=\sqrt{q_{v,d_i}}, $
and let $\mathbf{b}_v=(b_v(1),\ldots,b_v(p))^T$ and $D_v=diag(b_v(i))$.
Further let $\mathbf{a}$ denote the $p \times 1$ vector with entries
$
a_j = \sum_w \sqrt{{\bar{p}}_w}  (q_{w, d_j -1} - q_{w, d_j}).
$
This gives
$
\Sigma_0 = \sum_v D_v (I_p-\mathbf{b}_v\mathbf{b}_v^T) D_v + \mathbf{a} \mathbf{a}^T, \nn
$
where $I_p$ is the $p\times p$ identity matrix. For any matrix $A$, let $\rho_1(A) \leq \cdots \leq \rho_p(A)$ denote the eigenvalues of $A$ in increasing order.
By Weyl's Theorem (\cite{HornJohnson}, Theorem 4.3.1),
\beas
\rho_1(\Sigma_0) &\ge & \sum_v \rho_1(D_vB_vD_v) + \rho_1( \mathbf{a} \mathbf{a}^T) \ge \sum_v \rho_1(D_vB_vD_v).
\enas
Letting $B_v=I_p-\mathbf{b}_v\mathbf{b}_v^T$, it can be shown that the eigenvalues of $B_v$ are 1, with multiplicity $p-1$, and $\rho_1(B_v)=1-\mathbf{b}_v^T\mathbf{b}_v$ the least eigenvalue corresponding to the eigenvector $\mathbf{b}_v$. Now, using the Rayleigh-Ritz characterization of eigenvalues (\cite{HornJohnson}, Theorem 4.2.2),
\bea
\rho_1(D_vB_vD_v) &=& \min_{\mathbf{x}} \frac{ \mathbf{x}^TD_vB_vD_v\mathbf{x} }{ \mathbf{x}^T\mathbf{x} } = \min_{\mathbf{y}} \frac{ \mathbf{y}^TB_v\mathbf{y} }{ \mathbf{y}^TD_v^{-2}\mathbf{y} }\nn\\
&\geq& \frac{ \rho_1(B_v) }{ \rho_p(D_v^{-2}) } =\min_{i} q_{v,d_i}(1-\sum_{i}q_{v,d_i}).\nn
\ena
It therefore follows that
\bea
\|\Sigma_0^{-1/2}\| &\leq& \rho_p(\Sigma_0^{-1/2})= \frac{1}{\rho_1(\Sigma_0^{-1/2})} \\
&\leq & \Big[\sum_v \min_{i} q_{v,d_i}(1-\sum_{i}q_{v,d_i})\Big]^{-1/2} =: \tau.\nn
\ena
Finally we bound $\sum_{i,j} | \sigma_{ij} - \sigma_{ij}^0 | $. With \eqref{qdiff},
\beas
&& \left| \sigma_{ij} - \sigma_{ij}^0 \right|\\ 
&=& \Big | \sum_{v,w}  \Big\{ p_{wv}(1-p_{wv})(q_{v,d_i-1}^{(w)}-q_{v,d_i})(q_{w,d_j-1}^{(v)}-q_{w,d_j})  \nonumber \\
&& - \sqrt{{\bar{p}}_v} \sqrt{{\bar{p}}_w}(q_{v,d_i-1}-q_{v,d_i})(q_{w,d_j-1}-q_{w,d_j}) \Big\} \Big| \nn \\
&\le &
4 \sum_{v,w}    p_{wv}^2
 + \sum_{v; w \ne v }  | p_{wv} - \sqrt{{\bar{p}}_v} \sqrt{{\bar{p}}_w}| |(q_{v,d_i-1}-q_{v,d_i})(q_{w,d_j-1}-q_{w,d_j}) | \\
 &&+ \sum_v {\bar{p}}_v |  q_{v,d_i-1}-q_{v,d_i} | | q_{v,d_j-1}-q_{v,d_j}|  =: S.
\enas
Collecting the bounds and using that $q_{v, d_i} \le 1$ gives  the result.
\epf

\bpf [\textbf{Proof of Theorem \ref{PoissonThm}}] Recall that 
$$q_{v,i} = \mathbb{E}X_{(v,i)} = \mathbb{P}(D(v)=i).$$
Based on Theorem 10.B in Barbour, Holst $\&$ Janson (1992), we obtain that
\beas
&& d_{TV}\big(\mathcal{L}(\Xi_M), \mbox{Po}(\bm{\lambda}_M)\big)\\ 
&=&\sum_{v\in V}\sum_{i=M}^{n-1}q_{v,i}^2+\Bigg\{\sum_{v\in V}\sum_{i=M}^{n-1}q_{v,i}\sum_{j=M,j\neq i}^{n-1}\mathbb{E}|X_{(v,j)}-X_{(v,j)}^{(v,i)}|\nn\\
&&+\sum_{v\in V}\sum_{i=M}^{n-1}q_{v,i}\sum_{u\in V_v}\sum_{j=M}^{n-1}\mathbb{E}|X_{(u,j)}-X_{(u,j)}^{(v,i)}|\Bigg\}. 
\enas
Note that, the first summand in the bracket is for the case $u=v$ for $j\neq i$, and the second summand in the bracket is for $u\in V_v$ for all $j$, covering all $(u,j)$ except the case $(u,j)=(v,i)$. Since $X_{(v,j)}^{(v,i)}=0$ for $j\neq i$ and $\mathbb{E}X_{(v,j)}=q_{v,j}$, the first summand in the bracket is equal to $\Sigma_{v\in V}\Sigma_{i=M}^{n-1}q_{v,i}\Sigma_{j=M,j\neq i}^{n-1}q_{v,j}$. This, together with the summand outside the bracket, yields $\Sigma_{v\in V}[\Sigma_{i=M}^{n-1}q_{v,i}]^2$.

Let $\delta(u,v,i,j)=\mathbb{E}|X_{(u,j)}-X_{(u,j)}^{(v,i)}|$ for $u\in V_v$, then the construction \eqref{II} gives
\be \lefteqn{\delta(u,v,i,j)}\nn\\
&=&\sum_{d:d>i}q_{v,d}\sum_{\mathbf{x}_d\subset V_v}\mathbb{P}(N(v)=\mathbf{x}_d|D(v)=d)\mathbb{E}\Big[|X_{(u,j)}-X_{(u,j)}^{(v,i)}|\,\Big|\,N(v)=\mathbf{x}_d\Big]\nn\\
&&+\sum_{d:d<i} q_{v,d}\sum_{\mathbf{x}_d\subset V_v}\mathbb{P}(N(v)=\mathbf{x}_d|D(v)=d) \nn \\
&& \mathbb{E}\Big[|X_{(u,j)}-X_{(u,j)}^{(v,i)}|\,\Big|\,N(v)=\mathbf{x}_d\Big].\nn
\ee

Let $\delta_1(u,v,i,j)$ and $\delta_2(u,v,i,j)$ denote the above two terms respectively. We first calculate $\delta_1(u,v,i,j)$. Note that, if $D(v)=d>i$ and $N(v)=\mathbf{x}_d$, then by construction \ref{II}, $|X_{(u,j)}-X_{(u,j)}^{(v,i)}|=1$ if and only if vertex $u$ satisfies (1) $u\in\mathbf{x}_d$, (2) $D(u)=j+1$ or $j$ in $\mathscr{G}(n,\{p_{ij}\})$ and (3) $u\nsim v$ in $\mathscr{G}^{(v,i)}(n,\{p_{ij}\})$, namely $u$ is not adjacent to $v$ in the graph. Hence, $\delta_1(u,v,i,j)$ can be written as
{
\be \delta_1(u,v,i,j)
&=&\sum_{d:d>i}q_{v,d}\sum_{\mathbf{x}_d\subset V_v:u\in\mathbf{x}_d}\mathbb{P}(N(v)=\mathbf{x}_d\,|\,D(v)=d)\nn\\
&&\;\;\cdot\mathbb{E}\Big[\mathbbm{1}\{D(u)=j+1\;\mbox{or}\;j\}\mathbbm{1}^{(v,i)}\{u\nsim v\}\,\Big|\,N(v)=\mathbf{x}_d\Big].\nn
\ee
}
Conditioning on the event ``\,$D(u)=j+1 \mbox{ or }j$",
since the neighbourhood of vertex $v$ does not contain any relevant information for the degree of vertex $u$ once we know whether or not $u\sim v$, we have
\be
&&  \delta_1(u,v,i,j) \nn \\ 
&=&\sum_{d:d>i}q_{v,d} \sum_{\mathbf{x}_d\subset V_v:u\in\mathbf{x}_d}\mathbb{P}(N(v)=\mathbf{x}_d\,|\,D(v)=d)\label{d1}\\
&&\;\;\cdot\mathbb{P}(D(u)=j+1\;\mbox{or}\;j\,|\,X_{\{u,v\}}=1)\cdot
\mathbb{P}(X_{\{u,v\}}^{(v,i)}=0\,|\,N(v)=\mathbf{x}_d).\nn
\ee

If $d>i$ and $u\in\mathbf{x}_d$, then $\mathbb{P}(X_{\{u,v\}}^{(v,i)}=0\,|\,N(v)=\mathbf{x}_d)$ is the probability that the edge $\{u,v\}$ is  deleted in the construction \ref{II} and hence,  by \eqref{delp},
\be
\mathbb{P}(X_{\{u,v\}}^{(v,i)}=0\,|\,N(v)=\mathbf{x}_d)
&=& \sum_{\mathbf{x}_i\subset V_v:\mathbf{x}_i\subset\mathbf{x}_d,u\notin\mathbf{x}_i}\mathbb{P}(N^{(v,i)}(v)=\mathbf{x}_i\,|\,N(v)=\mathbf{x}_d) \nn \\
&=&\sum_{\mathbf{x}_i\subset V_v:\mathbf{x}_i\subset\mathbf{x}_d,u\notin\mathbf{x}_i}f^{+}(\mathbf{x}_i\,|\,\mathbf{x}_d).\label{SingleVertexDelete}
\ee
For (\ref{d1}), we then obtain
\be \delta_1(u,v,i,j)&=&\sum_{d:d>i}\;\sum_{\mathbf{x}_d\subset V_v:u\in\mathbf{x}_d}\;\sum_{\mathbf{x}_i\subset V_v:\mathbf{x}_i\subset\mathbf{x}_d,u\notin\mathbf{x}_i}
f^{+}(\mathbf{x}_i\,|\,\mathbf{x}_d)\nn\\
&&\cdot\mathbb{P}(N(v)=\mathbf{x}_d)\mathbb{P}(D(u)=j+1\;\mbox{or}\;j\,|\,X_{\{u,v\}}=1).\label{NewDelta1}
\ee

Next, we find $\delta_2(u,v,i,j)$ by a similar argument as for $\delta_1(u,v,i,j)$. If $D(v)=d<i$ and $N(v)=\mathbf{x}_d$, then by construction \ref{II}, $|X_{(u,j)}-X_{(u,j)}^{(v,i)}|=1$ if and only if vertex $u$ satisfies (1) $u\notin\mathbf{x}_d$, (2) $D(u)=j-1$ or $j$ in $\mathscr{G}(n,\{p_{ij}\})$ and (3) $u\sim v$ in $\mathscr{G}^{(v,i)}(n,\{p_{ij}\})$. Hence, following the same argument as for $\delta_1(u,v,i,j)$, we arrive at
\be
&& \delta_2(u,v,i,j)\nn
\\
&=&\sum_{d:d<i}\;\sum_{\mathbf{x}_d\subset V_v:u\notin\mathbf{x}_d}\;\sum_{\mathbf{x}_i\subset V_v:\mathbf{x}_i\supset\mathbf{x}_d,u\in\mathbf{x}_i}
f^{-}(\mathbf{x}_{i}\,|\,\mathbf{x}_{d})\nn\\
&&\cdot\mathbb{P}(N(v)=\mathbf{x}_d)\mathbb{P}(D(u)=j-1\;\mbox{or}\;j\,|\,X_{\{u,v\}}=0).\label{NewDelta2}
\ee

Now, the sum of at (\ref{NewDelta1})  (\ref{NewDelta2}) gives $\delta(u,v,i,j)$, and hence
$$ d_{TV}\big(\mathcal{L}(\Xi_M), \mbox{Po}(\bm{\lambda}_M)\big)\leq b_{M,1}+b_{M,2}+b_{M,3}, $$
where
\be
b_{M,1}=\sum_{v\in V}\bigg[\sum_{i=M}^{n-1}\mathbb{P}(D(v)=i)\bigg]^2,\nn
\ee
\be
b_{M,2}&=&\sum_{v\in V}\sum_{i=M}^{n-1}\sum_{u\in V_v}\sum_{j=M}^{n-1}\sum_{d:d>i}\;\sum_{\mathbf{x}_d\subset V_v:u\in\mathbf{x}_d}\;\sum_{\mathbf{x}_i\subset V_v:\mathbf{x}_i\subset\mathbf{x}_d,u\notin\mathbf{x}_i}
f^{+}(\mathbf{x}_i\,|\,\mathbf{x}_d)\nn\\
&&\cdot\mathbb{P}(D(v)=i)\mathbb{P}(N(v)=\mathbf{x}_d)\mathbb{P}(D(u)=j+1\;\mbox{or}\;j\,|\,X_{\{u,v\}}=1),\nn
\ee
and
\be
b_{M,3}&=&\sum_{v\in V}\sum_{i=M}^{n-1}\sum_{u\in V_v}\sum_{j=M}^{n-1}\sum_{d:d<i}\;\sum_{\mathbf{x}_d\subset V_v:u\notin\mathbf{x}_d}\;\sum_{\mathbf{x}_i\subset V_v:\mathbf{x}_i\supset\mathbf{x}_d,u\in\mathbf{x}_i}
f^{-}(\mathbf{x}_{i}\,|\,\mathbf{x}_{d})\nn\\
&&\cdot\mathbb{P}(D(v)=i)\mathbb{P}(N(v)=\mathbf{x}_d)\mathbb{P}(D(u)=j-1\;\mbox{or}\;j\,|\,X_{\{u,v\}}=0).\nn
\ee


In particular, removing the constraints of vertex $u$ on choosing $\mathbf{x}_d$ and $\mathbf{x}_i$, and using Lemma \ref{isaprob},
we obtain that
\be
b_{M,2}&\leq & 2\sum_{v\in V}\sum_{u\in V_v}\sum_{i=M}^{n-1}
\mathbb{P}(D(v)=i)\,\mathbb{P}(D(v)>i)\,\mathbb{P}(D(u)\geq M\,|\,X_{\{u,v\}}=1) \nn\\
&\leq& 2\sum_{v\in V}\sum_{u\in V_v}
\mathbb{P}(D(v)\geq M)\,\mathbb{P}(D(u)\geq M\,|\,X_{\{u,v\}}=1) \nn
\ee
and similarly
$
 b_{M,3}
\leq  2\sum_{v\in V}\sum_{u\in V_v}
\mathbb{P}(D(v)\geq M)\,\mathbb{P}(D(u)\geq M-1\,|\,X_{\{u,v\}}=0).
$
Now, with both $\mathbb{P}(D(u)\geq M\,|\,X_{\{u,v\}}=1)$ and $\mathbb{P}(D(u)\geq M-1\,|\,X_{\{u,v\}}=0)$ equal to $\mathbb{P}(D^{(v)}(u)\geq M-1)$, the assertion follows.

\epf

\textbf{Acknowledgements.} {GR was supported in part by EPSRC EP/D003/05/1 Amorphous Computing Grant, and by BBSRC and EPSRC through OCISB.}


\begin{thebibliography}{99}
\footnotesize

\bibitem{barabasialbert} {\sc Barabasi, A.~L.and Albert, R.} (1999). Emergence of scaling in random network. {\it Science} 286, 509--512.

\bibitem{bhj} {\sc Barbour, A.~D., Holst, L. and Janson, S. } (1992).
{\sl Poisson Approximation.} Oxford University Press. 


\bibitem{larryjay} {\sc Bartroff, J. and  Goldstein, L.} (2009). A Berry-Esseen bound with applications to the number of multinomial cells of given occupancy and the number of graph vertices of given degree. Preprint.


\bibitem{bollobasetal} {\sc Bollob\'as, B., Janson, S. and Riordan, O.} (2007). The phase transition in inhomogeneous random graphs.
{\it Random Structures and Algorithms} {\bf 31}, 3--122.

\bibitem{bbb5} {\sc Bollob\'{a}s}, B. (2001). \emph{Random Graphs}. Cambridge University Press, Cambridge.


\bibitem{daudinetal} {\sc  Daudin, J-J., Picard, F. and Robin, S.} (2008).
A mixture model for random graphs. {\it Statistics and Computing} {\bf 18}, 173--183.

\bibitem{dormenbook} {\sc Dorogovtsev, S.~N. and Mendes, J.~F.~F.}  (2003).
{\it Evolution of Networks:
From Biological Nets to the Internet and WWW\/}.
Oxford University Press, Oxford.



\bibitem{er} {\sc Erd\"{o}s, P. and R\'{e}nyi, A.} (1959). On random graphs. \emph{Publicationes Mathematicae} \textbf{6}, 290-297.


\bibitem{GoldsteinRinott1996}  {\sc Goldstein, L. and Rinott, Y.} (1996).
Multivariate normal approximations by Stein's method and size bias couplings.
\emph{J. Appl. Probab.} { \bf 33},  1--17.

 \bibitem{HornJohnson}  {\sc Horn, R.~A. and Johnson,  C.~R.} (1985). \textit{Matrix Analysis}. Cambridge University Press.

\bibitem{lin}  {\sc Lin, K.} (2008). \textit{Motif counts, clustering coefficients, and vertex degrees in models of random networks.}  DPhil Dissertation, Department of Statistics, University of Oxford.


\bibitem{a} {\sc McKay, B.~D. and Wormald, N.~C.} (1997). The degree sequence of a random graph. I. The models. \emph{Random Structures and Algorithms} \textbf{11}, 97-117.


\bibitem{nmw} {\sc Newman,  M.~E.~J., Moore, C. and Watts, D.~J.} (2000). Mean-field solution of the small-world network model. \emph{Physical Review Letters} \textbf{84}, 3201-3204.

\bibitem{nowickisnijders} {\sc Nowicki, K. and Snijders, T.~A.~B.} (2001) Estimation and prediction for stochastic blockstructures. \emph{JASA} { \bf{96}}, 1077-1087.


\bibitem{Rinott1996} {\sc Rinott, Y. and Rotar, V.} (1996).
 A multivariate {CLT} for local dependence with {$n\sp {-1/2}\log n$}
  rate and applications to multivariate graph related statistics. \emph{J. Multivariate Analysis}  
\textbf{ 56}, \penalty0 333--350.

\bibitem{solowetal}{\sc Solow, A.~R., Costello,  C.~J. and Ward, M.} (2003). Testing the power law model for discrete size data. {\it Ann. Nat.} 162, 685--689.

\bibitem{Stein_orig}  {\sc Stein, C.} (1972). A bound for the error in the normal approximation to the distribution of a sum of dependent random variables. In:
    \textit{Proceedings of the Sixth Berkeley Symposium on Mathematical Statistics and Probability, Vol. II: Probability theory}, 583-602. Univ. California Press,
    Berkeley, CA.

\bibitem{Stein_book}  {\sc Stein, C.} (1986). \textit{Approximate computation of expectations.} Institute of Mathematical Statistics Lecture Notes - Monograph Series,
    \textbf{7}. Institute of Mathematical Statistics, Hayward, CA.


\bibitem{stumpfetal}{\sc Stumpf, M.~P.~H., Wiuf,  C. and May, R.~M.} (2005). Subnets of scale-free networks are not scale-free: Sampling properties of networks. \emph{PNAS} 102, no.12, 4221-4224.


\end{thebibliography}
\end{document}